\documentclass[11pt]{article}
\usepackage{mathrsfs,amsmath,amsfonts,amssymb}
\usepackage{graphicx}
\usepackage{color}
\usepackage{array}
\usepackage{lscape}
\usepackage{hyperref}

\usepackage{amsmath,amsfonts,amssymb}
\usepackage{multirow}

\renewcommand{\vec}[1]{\boldsymbol{#1}}

\newtheorem{remark}{Remark}[section]
\newtheorem{theorem}[remark]{Theorem}

\newtheorem{corollary}[remark]{Corollary}

\newcommand{\half}{\ensuremath{\frac{1}{2}}}

\newcommand{\ignore}[1]{}

\newcommand{\Fig}[1]{Figure \ref{fig:#1}}
\usepackage{siunitx}
\usepackage{booktabs}

\title{{The Cartesian Grid Active Flux Method with Adaptive Mesh Refinement}
\footnote{This work was funded by the Deutsche Forschungsgemeinschaft (DFG, German Research Foundation) - 325695158. Furthermore, research was supported by a subcontract with Boise State University, Boise, ID, under DARPA Cooperative Agreement HR00112120003 with Embry-Riddle Aeronautical University, Daytona Beach, FL, USA. This work is approved for public release; distribution is unlimited. The content of the information does not necessarily reflect the position or the policy of the Government, and no official endorsement should be inferred.}}
\author{Donna Calhoun\thanks{Department of Mathematics, Boise State
    University, Boise, ID, USA; {\tt  donnacalhoun@boisestate.edu}}
    \and Erik Chudzik\thanks{Institute of
    Mathematics, Heinrich-Heine-University D\"usseldorf,
    Germany} \thanks{{\tt Erik.Chudzik@hhu.de}}
   \and Christiane Helzel\footnotemark[3] \thanks{{\tt
       Christiane.Helzel@hhu.de} (corresponding author)}
 }

\date{}

\begin{document}
\maketitle

\begin{abstract}
We present the first implementation of the Active Flux method on adaptively refined Cartesian grids.

The Active Flux method is a third order accurate finite volume method for hyperbolic conservation laws, which is based on the use of point values as well as cell average values of the conserved quantities. The resulting method has a compact stencil in space and time and good stability properties.

The method is implemented as a new solver in ForestClaw, a software for parallel adaptive mesh refinement of patch-based solvers.  On each Cartesian grid  patch the single grid Active Flux method can be applied. The exchange of data between grid patches is organised via ghost cells. The local stencil in space and time and the availability of the point values that are used for the reconstruction, leads to an efficient implementation. The resulting method is third order accurate, conservative and allows the use of subcycling in time.       
\end{abstract}




\section{Introduction}
The Active Flux method is a finite volume method for hyperbolic conservation laws that was previously introduced by  Eymann, Roe and coauthors \cite{article:ER2011b,article:ER2011a,article:ER2013,article:Roe2017,article:RMD2018}. In its original form the method is third order accurate. This is achieved by using a continuous, piecewise quadratic reconstruction and a sufficiently accurate quadrature rule to compute numerical fluxes. The quadrature method, i.e.\ typically Simpson's rule,  requires point values of the conserved quantities  at grid cell interfaces at the current time as well as at later time levels. These point values together with the cell average value are also used to compute the reconstruction. While classical finite volume methods only use cell average values of the conserved quantities as degrees of freedom, the Active Flux method involves both point values and cell average values degrees of freedom. This adds flexibility to the numerical method.  For linear advection  and the acoustic equations in one, two or three spatial dimensions, these point values can be updated using exact evolution formulas, making the method truly multidimensional. 

Originally, Roe and Ey\-mann \cite{article:ER2013} used un\-structured tri\-angular grids for their two-dimensional active flux method. In \cite{article:BHKR2018} and \cite{article:HKS2019}, two-dimen\-sio\-nal Cartesian grid versions of the active flux method were introduced. The Cartesian grid method will also be used in this paper and briefly reviewed in the next section. 

While the order of convergence of a numerical scheme is a property that can be shown in the limit when the mesh width and the time step goes to zero, for practical computations it is desirable to obtain accurate results on relatively coarse grids. Roe \cite{article:Roe2021} argues that the accurate approximation on coarse grids is strongly influenced by the computational stencil and that exact evolution operators perform well in this respect. Barsukow  showed that the two-dimensional Cartesian grid Active Flux method  for the acoustic equations with an exact evolution operator as described in  \cite{article:Barsukow2020} is stationary preserving. This means that the numerical scheme does not add dissipation to discrete representations of all stationary states of the acoustic equations. As a consequence such states can be computed with very high accuracy on coarse grids, while other methods would require a much higher resolution.

In earlier related work,  Luk\'a\v{c}ov\'a-Medvid'ov\'a et al.  \cite{article:LMW2000,article:LSW2002}, use exact multidimensional evolution operators as building blocks of finite volume methods of various order. Those methods differ from the Active Flux method mainly in the choice of the degrees of freedom.

Although the Active Flux method performs well on coarse grids, for practical applications it might still be desirable to vary the size of the grid cells adaptively in order to allow a higher resolution in parts of the computational domain as needed. A possible application that could benefit from local refinement is the propagation of a high frequency acoustic wave. 

In this paper we show how the Active Flux method can be applied on adaptively refined Cartesian grids. The local stencil of the method allows an efficient transfer of data between the different grid patches. Our numerical results confirm third order accuracy of the resulting method.The method is implemented as a new solver in ForestClaw  \cite{article:BCMT2014,misc:ForestClaw,ca-bu:2017}.

This paper is organised as follows. In Section \ref{sec:AF} we briefly review the Active Flux method for two-dimensional Cartesian grids. Section \ref{sec:3} describes the extension to adaptively refined grids. In Section \ref{sec:4} we introduce new Active Flux methods for advective transport and illustrate the performance of the adaptively refined Active Flux method for a variety of test problems. 

\section{The Cartesian grid Active Flux Method}
\label{sec:AF}
In this section we provide a brief review of the Active Flux method on a single two-dimensional Cartesian grid. More details can be found in  Barsukow et al.\ \cite{article:BHKR2018} and Helzel et al.\ \cite{article:HKS2019}.

\begin{figure}[htb]
\centerline{\includegraphics[width=0.3\textwidth]{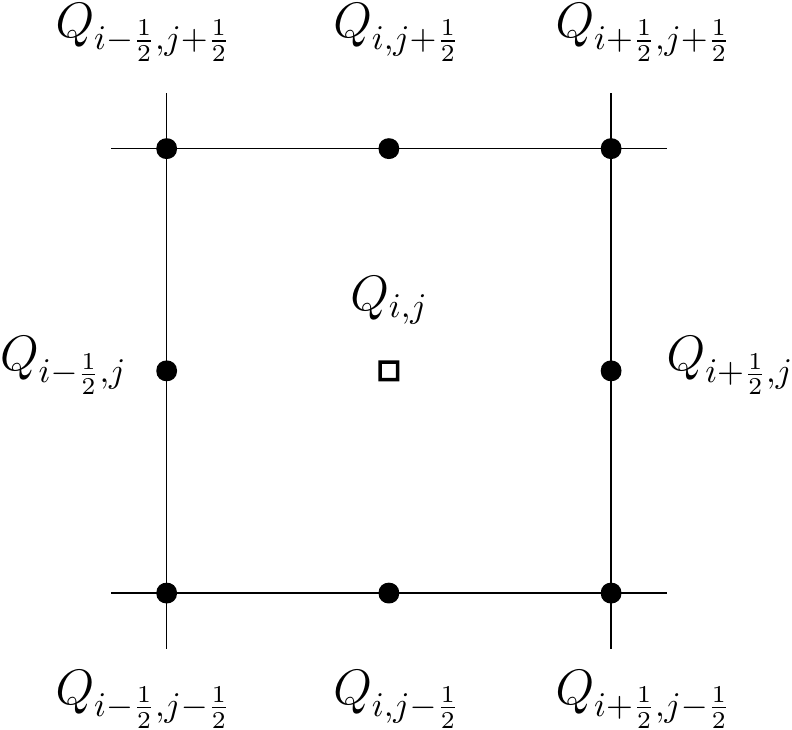}} 
\caption{\label{fig:dof} Degrees of freedom used to reconstruct the piecewise quadratic function in grid cell $(i,j)$. Point values are marked by a dot and the cell average is marked by a square.}
\end{figure}

We consider hyperbolic conservation laws in divergence form
\begin{equation*}
\partial_t q + \partial_x f(q) + \partial_y g(q) = 0,
\end{equation*}
where $q:\mathbb{R}^2 \times \mathbb{R}^+ \rightarrow \mathbb{R}^s$ is a vector of conserved quantities and $f, g : \mathbb{R}^s \rightarrow \mathbb{R}^s$ are vector valued flux functions. On a single patch we use a two-dimensional Cartesian grid with equidistant mesh sizes $\Delta x$ and $\Delta y$. The grid cell $(i,j)$ is described by $[ x_{i-\frac{1}{2}},x_{i+\frac{1}{2}}] \times [ y_{j-\frac{1}{2}}, y_{j+\frac{1}{2}}] \subset \mathbb{R}^2, \, i,j, \in \mathbb{Z}.$ As a finite volume method, the Active Flux method computes  cell averaged values of the conserved quantities via an update of the form
\begin{equation}
\label{eqn:2dfv}
Q_{i,j}^{n+1} = Q_{i,j}^n - \frac{\Delta t}{\Delta x} \left(
  F_{i+\frac{1}{2},j} - F_{i-\frac{1}{2},j} \right) - \frac{\Delta
  t}{\Delta y} \left( G_{i,j+\frac{1}{2}} - G_{i,j-\frac{1}{2}} \right),
\end{equation}
where $Q_{i,j}^n$ is an approximation of the cell average values of the conserved quantities in grid cell $(i,j)$ at time $t_n$ and $F_{i\pm \frac{1}{2},j}$, $G_{i,j\pm \frac{1}{2}}$ are numerical fluxes at the grid cell interfaces given by
\begin{equation*}
\begin{split}
F_{i+\frac{1}{2},j} & \approx \frac{1}{\Delta t \Delta y} \int_{t_n}^{t_{n+1}}
\int_{y_{j-\frac{1}{2}}}^{y_{j+\frac{1}{2}}}
f(q(x_{i+\frac{1}{2}},y,t)) dy \, dt\\
G_{i,j+\frac{1}{2}} & \approx \frac{1}{\Delta t \Delta x}
\int_{t_n}^{t_{n+1}} \int_{x_{i-\frac{1}{2}}}^{x_{i+\frac{1}{2}}}
g(q(x,y_{j+\frac{1}{2}},t)) dx \, dt. 
\end{split}
\end{equation*}  
As suggested by Eymann and Roe, we use Simpson's rule to compute the numerical fluxes.  For fluxes $F_{i + \frac{1}{2},j}$, this leads to the formula
\begin{equation}\label{eqn:Simpson-2d}
  \begin{split}
F_{i+\frac{1}{2},j} & := \frac{1}{36} \Big(
  f(Q_{i+\frac{1}{2},j-\frac{1}{2}}^n) + 4 f(Q_{i+\frac{1}{2},j}^n) +
  f(Q_{i+\frac{1}{2},j+\frac{1}{2}}^n)  \\
  & \quad  + 4
    f(Q_{i+\frac{1}{2},j-\frac{1}{2}}^{n+\frac{1}{2}})  + 16   f(Q_{i+\frac{1}{2},j}^{n+\frac{1}{2}}) + 4
    f(Q_{i+\frac{1}{2},j+\frac{1}{2}}^{n+\frac{1}{2}}) \\ 
    & \quad  + f(Q_{i+\frac{1}{2},j-\frac{1}{2}}^{n+1}) + 4
     f(Q_{i+\frac{1}{2},j}^{n+1}) +
     f(Q_{i+\frac{1}{2},j+\frac{1}{2}}^{n+1}) \Big)
  \end{split}  
\end{equation}
We use an analogous formula for flux $G_{i,j+\frac{1}{2}}$. The $Q$ values in the right hand side of (\ref{eqn:Simpson-2d}) are approximations to point values of the conserved quantities at the grid cell interface  at times $t_n$, $t_{n+\half}$ and $t_{n+1}$.  To compute these point values, we assume that at time $t_n$ the average values $Q_{i,j}^n$ and point values at cell corners $Q_{i\pm \frac{1}{2}, j \pm \frac{1}{2}}^n$ and edge midpoints $Q_{i\pm \frac{1}{2},j}^n$ and $Q_{i,j\pm \frac{1}{2}}^n$ are known.  The location of these point values are shown in Figure \ref{fig:dof}.
Based on these known values, a two-dimensional quadratic polynomial can be reconstructed in each grid cell. For simplicity, all Cartesian grid cells are mapped to a reference cell $[-1,1] \times [-1,1]$. In each cell, the reconstruction has the form
\begin{equation}
q^{rec}(\xi,\eta) = c_{00} + c_{10} \xi + c_{01} \eta + c_{20} \xi^2 +
c_{11} \xi \eta + c_{02} \eta^2 + c_{21} \xi^2 \eta + c_{12} \xi
\eta^2 + c_{22} \xi^2 \eta^2
\label{eqn:poly}
\end{equation}  
with $c_{ij} \in \mathbb{R}$, $i,j = 0, 1, 2$ and $(\xi,\eta) \in [-1,1] \times [-1,1]$. The precise form of the reconstructed function is obtained by interpolating the eight known point values along the grid cell boundary and by requiring that the average of the reconstructed function agrees in each grid cell with the known cell average. This can conveniently be obtained by using appropriate basis functions as introduced in \cite{article:HKS2019}. 

We denote the (at time $t_n$) reconstructed continuous,  piecewise quadratic function by $q^n(x,y)$. The quadratic polynomial that was reconstructed in grid cell $(i,j)$ is denoted by $q_{ij}^n(x,y)$.

For special linear partial differential equations, the required point values of the conserved quantities  can be calculated  using an exact evolution formula applied to the piecewise quadratic data.  This is in particular the case for linear advection and acoustics. The use of the exact evolution formula leads to a truly multidimensional numerical method. While the point values at time $t_{n+\frac{1}{2}}$ are only used in order to compute the numerical fluxes, the point values at time $t_{n+1}$ are also used to compute the continuous, piecewise quadratic reconstruction $q^{n+1}(x,y)$, which
provides the initial data for the next time step.

Since the point values are located along the grid cell boundaries, they are used for the reconstruction in two or four grid cells.  Thus, the Active Flux method for two-dimensional Cartesian grids has four degrees of freedom  per grid cell: the cell average plus three point values along the grid cell boundary. In our implementation, we associate point values $Q_{i-\frac{1}{2},j}$, $Q_{i-\frac{1}{2},j-\frac{1}{2}}$, $Q_{i,j-\frac{1}{2}}$  and the average value $Q_{i,j}$ with grid cell $(i,j)$.

\subsection{Advection equation}
\label{sec:2_1}
For the advection equation
\begin{equation}
\partial_t q + a \partial_x q + b \partial_y q = 0, \label{eqn:advection2d}
\end{equation}
with $q:\mathbb{R}^2 \times \mathbb{R}^+ \rightarrow \mathbb{R}$, $a,
b \in \mathbb{R}$,
the exact evolution formula has the form
$$
q(x,y,t+\tau) = q(x-a\tau,y-b \tau, t).  
$$
Thus, the required point values at time $t_{n+\frac{1}{2}}$ and
$t_{n+1}$ can easily be computed by
evaluating the continuous, piecewise quadratic function $q^n$ at the
appropriate upwind points, for example
\begin{equation}\label{eqn:advection-update}
\begin{split}
 Q_{i+\frac{1}{2},j- \frac{1}{2}}^{n+\frac{1}{2}} & =
q^n\left( x_{i+\frac{1}{2}}-a \frac{\Delta t}{2},y_{j-\frac{1}{2}}-b
  \frac{\Delta t}{2} \right) \\
Q_{i+\frac{1}{2},j- \frac{1}{2}}^{n+1} & =
q^n \left( x_{i+\frac{1}{2}}-a \Delta t,y_{j-\frac{1}{2}}-b
\Delta t \right).
\end{split} 
\end{equation}  

In the numerical method, we want to restrict the time steps, so that waves propagate at most one grid cell per time step. This can be expressed in the form
\begin{equation}\label{eqn:cfl-adv}
\mbox{CFL}:=\mbox{max} \left( \frac{|a| \Delta t}{\Delta x}, \frac{|b| \Delta
    t}{\Delta y} \right) \le 1. 
\end{equation}  
We have shown in \cite{article:CHK2021} that linear stability requires a more severe time step restriction. In practical computations we therefore use time steps which satisfy $\mbox{CFL} \le 0.75$.

In Section \ref{sec:Numerics-advection} the approach is extended to advective transport in a spatially and temporally varying velocity field.

\subsection{Burgers' equation}\label{sec:2-Burgers}
For the two-dimensional Burgers' equation
\begin{equation}\label{eqn:burgers2d}
q_t +  \left( \frac{1}{2} q^2 \right)_x + \left( \frac{1}{2} q^2
\right)_y = 0
\end{equation}
with $q : \mathbb{R}^2 \times \mathbb{R}^+ \rightarrow \mathbb{R}$, we do not have an exact evolution formula. Instead we  approximate the point values using an approximative approach described in \cite{article:CHK2021}. 

For smooth solutions, equation (\ref{eqn:burgers2d}) can equivalently be written in the advective form
\begin{equation}
q_t + q q_x + q q_y = 0,  
\end{equation}
which suggests the implicitly defined evolution formula 
\begin{equation} \label{eqn:burgers-sol}
    q(x,y,t)  = q(x-q(x,y,t)t,y-q(x,y,t)t,0). 
\end{equation}  
Starting with an initial guess 
\begin{equation*}
\left( Q_{i-\frac{1}{2},j}^{n+\frac{1}{2}}\right)^0, \quad \left( Q_{i-\frac{1}{2},j}^{n+1}\right)^0,
\end{equation*}
we iteratively compute
\begin{equation}\label{eqn:Burgers-iterativ}
\begin{split}  
\left( Q_{i-\frac{1}{2},j}^{n+\frac{1}{2}} \right)^\ell & = q^n \left( x_{i-\frac{1}{2}} -
  \left( Q_{i-\frac{1}{2},j}^{n+\frac{1}{2}} \right)^{\ell-1}
  \frac{\Delta t}{2}, y_j - \left( Q_{i-\frac{1}{2},j}^{n+\frac{1}{2}}
  \right)^{\ell-1} \frac{\Delta t}{2} \right) \\
\left( Q_{i-\frac{1}{2},j}^{n+1} \right)^\ell & = q^n \left( x_{i-\frac{1}{2}} -
  \left( Q_{i-\frac{1}{2},j}^{n+1} \right)^{\ell-1}
  \Delta t, y_j - \left( Q_{i-\frac{1}{2},j}^{n+1}
  \right)^{\ell-1} \Delta t \right) \quad \ell = 1, 2, \ldots.
\end{split}
\end{equation}
Wave speeds at other positions are computed analogously. Each iteration improves the accuracy by one order. We start the iteration with a first order accurate approximation.  The piecewise quadratic reconstruction limits the achievable accuracy to third order. Thus, a third order accurate approximation can be achieved after two iterations. 

An obvious initial guess would be to use the point values at the respective location, i.e.\
\begin{equation*}
\left( Q_{i-\frac{1}{2},j}^{n+\frac{1}{2}} \right)^0 = \left(
  Q_{i-\frac{1}{2},j}^{n+1} \right)^0
= Q_{i-\frac{1}{2},j}^n
\end{equation*}
and analogously for all other point values along the grid cell boundaries. These wave speeds are third order accurate in space and first order accurate in time. However, this choice suffers from an instability, if the characteristic speed changes sign as explained in detail in \cite{article:HKS2019}. The instability can be mitigated if data from all adjacent grid cells are used to compute the initial guess. This leads to a stronger coupling of wave speeds and cell average values. Here, we compute the initial guess of the wave speed based on the neighboring cell average values, i.e.\ we use
\begin{equation*}
  \begin{split}
    \left( Q_{i-\frac{1}{2},j}^{n+\frac{1}{2}} \right)^0 & =
    \left( Q_{i-\frac{1}{2},j}^{n+1} \right)^0= \frac{1}{2}
\left(Q_{i-1,j}^n+Q_{i,j}^n \right)\\
\left( Q_{i-\frac{1}{2},j-\frac{1}{2}}^{n+\frac{1}{2}} \right)^0 & =
\left( Q_{i-\frac{1}{2},j-\frac{1}{2}}^{n+1} \right)^0 =
\frac{1}{4}
\left( Q_{i-1,j}^n+Q_{i,j}^n+Q_{i-1,j-1}^n+Q_{i,j-1}^n \right) \\
\left( Q_{i,j-\frac{1}{2}}^{n+\frac{1}{2}} \right)^0 & =
\left( Q_{i,j-\frac{1}{2}}^{n+1} \right)^0= \frac{1}{2}
\left(
  Q_{i,j}^n + Q_{i,j-1}^n \right).
\end{split}
\end{equation*}
The slight increase of the stencil, as introduced by this initial guess of the wave speeds, leads to a stable approximation.

\subsection{Acoustics}
\label{sec:2_2}
The acoustic equations are given by
\begin{equation}\label{eqn:acoustics}
  \begin{split}
  \partial_t p + c \nabla \cdot \vec{u} & = 0 \\
  \partial_t \vec{u} + c \nabla p & = 0,
  \end{split}
\end{equation}
where $\vec{u}:\mathbb{R}^2 \times \mathbb{R}^+ \rightarrow \mathbb{R}^2$ is the velocity vector, $p:\mathbb{R}^2 \times \mathbb{R}^+ \rightarrow \mathbb{R}$ is the pressure and $c \in \mathbb{R}^+$ is the speed of sound.

The evolution formula for the two-dimensional acoustic equations used in this paper can be found in \cite{article:ER2013}. It is based on the observation that (\ref{eqn:acoustics}) can be rewritten as
\begin{equation}
  \begin{split}
    \partial_{tt} p - c^2 \Delta p & = 0 \\
    \partial_{tt} \vec{u} - c^2 \Delta \vec{u} & = c^2 \nabla \times \vec{w},
  \end{split}
\end{equation}
where $\vec{w} = \nabla \times \vec{u}$ is the vorticity and $\triangle$ is the Laplacian operator. In the two-dimensional case, considered here, $\vec{w} = (0,0,v_x-u_y)^T$. Thus, in a flow with constant vorticity both pressure and velocity satisfy a wave equation. Furthermore, it is easy to verify that the vorticity is stationary, i.e.\ in the two-dimensional case the relation $\partial_t (v_x - u_y) = 0$ holds. 

Assuming constant vorticity, Eymann and Roe \cite{article:ER2013} derived the evolution formulas
\begin{equation}\label{eqn:ERFormula}
  \begin{split}
p(t) & = M_R \{ p \} + R \left( \partial_R M_R \{ p \} - M_R \{ \nabla
  \cdot \vec{u} \} \right) \\
\vec{u}(t) & = M_R \{ \vec{u} \} + R \left( \partial_R M_R \{ \vec{u}
  \} - M_R \{ \nabla p \} \right),
  \end{split}
\end{equation}
where $R = c \cdot t$ and $M_R\{ f\}$ is the spherical mean. The values of pressure and velocity at the right hand side of the evolution equation  (\ref{eqn:ERFormula}) are given initial values at time $t=0$.  For a scalar function $f:\mathbb{R}^2 \rightarrow \mathbb{R}$, the spherical mean over a disc with radius $R$,  centred at $(x,y)$ is defined by  
\begin{equation}\label{eqn:sphe-mean}
M_R \{ f \}(x,y) := \frac{1}{2 \pi R} \int_0^{2 \pi} \int_0^R f(x+s \cos
\phi, y+s \sin \phi) \frac{s}{\sqrt{R^2 - s^2}} ds d\phi.
\end{equation}
In the vector valued case, the formula is applied component wise. The solution formula can be evaluated exactly, if during each time step the previous values of pressure and velocity are replaced by the corresponding components of the reconstructed continuous, piecewise quadratic function $q^n$.

Each time step of the explicit Active Flux method is restricted so that the circle around the edge midpoint over which the integration takes place remains inside the two adjacent grid cells.  This condition will be met if  
\begin{equation}\label{eqn:cfl-acoustics}
\max \left( \frac{c \Delta t}{\Delta x}, \frac{c \Delta t}{\Delta y}
\right) \le \frac{1}{2}.
\end{equation}
In \cite{article:CHK2021} we showed that this necessary condition is sufficient for linear stability on a regular Cartesian grid.

\section{Adaptive Mesh Refinement for the Active Flux Method}
\label{sec:3}
We now describe the implementation of the adaptive Active Flux method as a new solver in ForestClaw \cite{article:BCMT2014}, a software for parallel adaptive mesh refinement based on a quadtree approach. In ForestClaw, Cartesian grid patches occupy quadrants in a quadtree, or multi-block forest of quadtrees. ForestClaw was developed by Calhoun and Burstedde based on the p4est software \cite{article:BWG2011}.

\begin{figure}[htb!]
    \centering
    \includegraphics[width=0.5\textwidth]{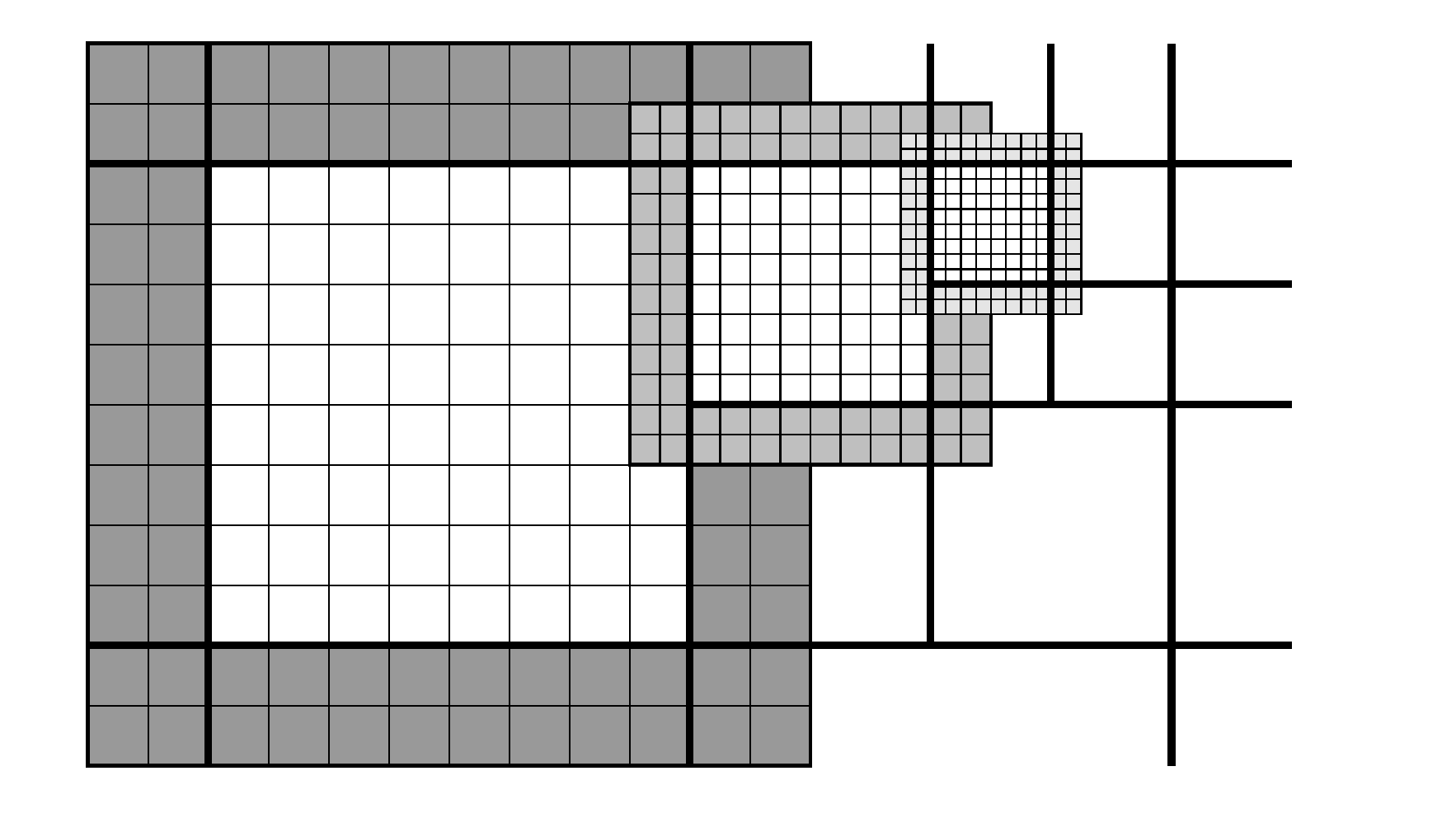}
    \caption{\label{fig:ForestClawGrid}
      Three refinement levels of quadrants of a ForestClaw mesh are shown. 
      Each quadrant is occupied by a logically Cartesian grid of fixed resolution (e.g. $8 \times 8$ in the figure).  Each grid also contains a layer of ghost cells (shaded region) which are used to facilitate the exchange of data between adjacent grid patches and between local and remote processors.}
  \end{figure}

Mesh refinement is realised by a bisection of grid patches so that a quadrant of resolution level $\ell$ is replaced by four quadrants of resolution level $\ell +1$. A patch of level zero would correspond to a single Cartesian grid discretizing a single, square domain. The number of grid cells on a single quadrant is constant for all levels, resulting in a 2:1 refinement ratio between resolution levels. Typically $8 \times 8$, $16 \times 16$ or $32 \times 32$ grid cells are used on a single patch, which allows  a flexible change of the resolution. Furthermore, single grid patches can efficiently be handled by separate processors in a parallel computation. Figure \ref{fig:ForestClawGrid} shows a typical situation. 

A ForestClaw mesh inherits many properties from the underlying p4est mesh \cite{article:BWG2011}. In particular, a ForestClaw mesh is {\em well-balanced} so that adjacent quadrants never differ by more than one level. This is a necessary condition for the definition of data transfer between different patches, which is organized using ghost cells and described in more detail below. Another key feature of the ForestClaw mesh is that it is dynamically adapted to follow solution features of interest.  At each time step, we apply coarsening and refinement criteria to the solution on each quadrant.  If the coarsening criteria is satisfied by the solution in each quadrant in a family of four quadrants,  the four quadrants will be replaced by a single quadrant, and the solution  will be averaged from the finer grids to the new coarser solution.  Otherwise, if the solution on a quadrant satisfies the refinement criteria, the quadrant will be subdivided into four quadrants, and the solution will be interpolated from the coarse parent to the new finer resolution solution on each of the child quadrants.  The details of the interpolation and averaging are provided next. 

\subsection{Spatial transfer of grid cell information}
Solution data in a composite ForestClaw mesh needs to be communicated between adjacent grids sharing quadrant boundaries, and when dynamically coarsening and refining the mesh. The Active Flux method allows a very efficient transfer of both pointwise and cell-average information. In the following, a "grid patch" is both the p4est quadrant and the solution data in the quadrant.  The following three situations need to be considered:
\begin{enumerate}
\item A \textbf{transfer from a fine grid to a coarse grid} is needed
  if four grid patches at level $\ell +1$ are coarsened to a single
  patch at level $\ell$. The same approach is used for the
  computation of ghost cell values on a patch of level $\ell$ from a
  neighbouring patch of level $\ell +1$.
\item A \textbf{transfer from a coarse grid to a fine grid} is needed
  if a patch of level $\ell$ is marked for refinement and four patches
  at level $\ell +1$ need to be reconstructed. The same approach is
  used for the computation of ghost cell values for a patch at level
  $\ell +1$ from a neighbouring patch of level $\ell$.
\item For neighbouring grid patches of the same level the ghost cell
  information is simply copied from the neighbouring grid patches.  
\end{enumerate}
We will now discuss the first two approaches in more detail.

\subsubsection*{Transfer from fine to coarse grids}
\begin{figure}[htb]
	\centerline{\includegraphics[width=0.8\textwidth]{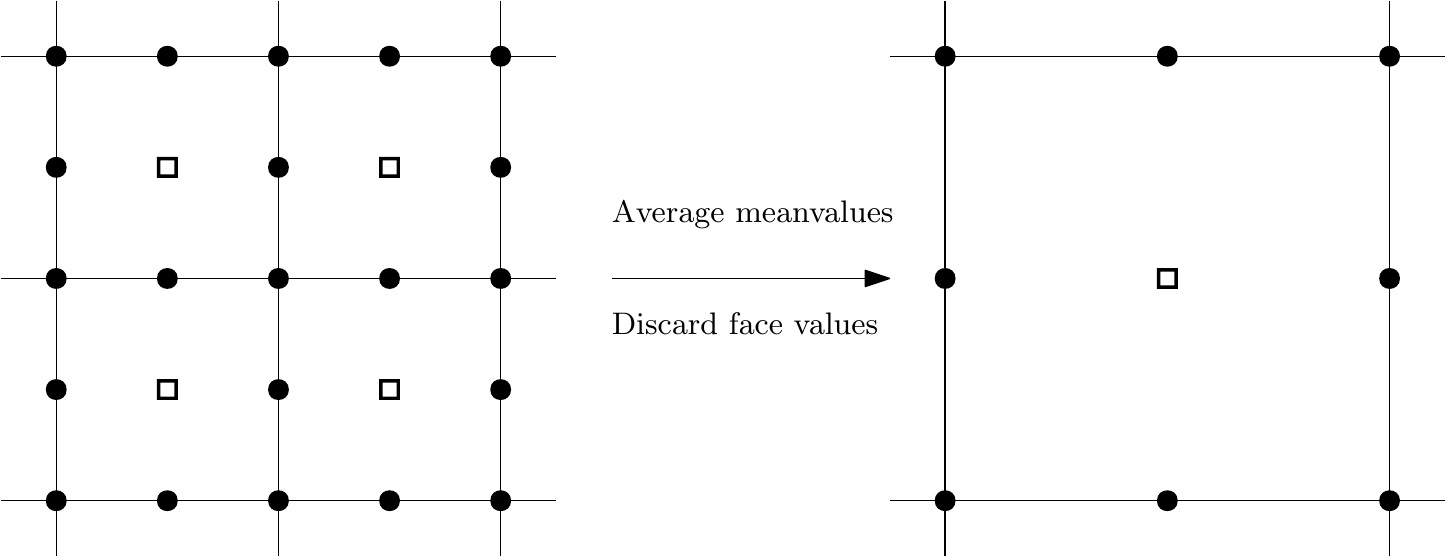}} 
	\caption{\label{fig:transfer1}Illustration of coarsening: The cell average on the
          coarse grid is the average of the four average values on the
          fine grid. The point values are copied from  vertices of
          the fine grid cells.}
\end{figure}

The degrees of freedom of a coarse grid cell are computed from the degrees of freedom of four grid cells on the finer level as illustrated in Figure \ref{fig:transfer1}. The cell average of the coarse cell is the average of the four cell average values on the fine grid. The point values are copied from the point values at the vertices of the fine grid cells.  

\subsubsection*{Transfer from coarse to fine grids}
\begin{figure}[htb]
	\centerline{\includegraphics[width=0.8\textwidth]{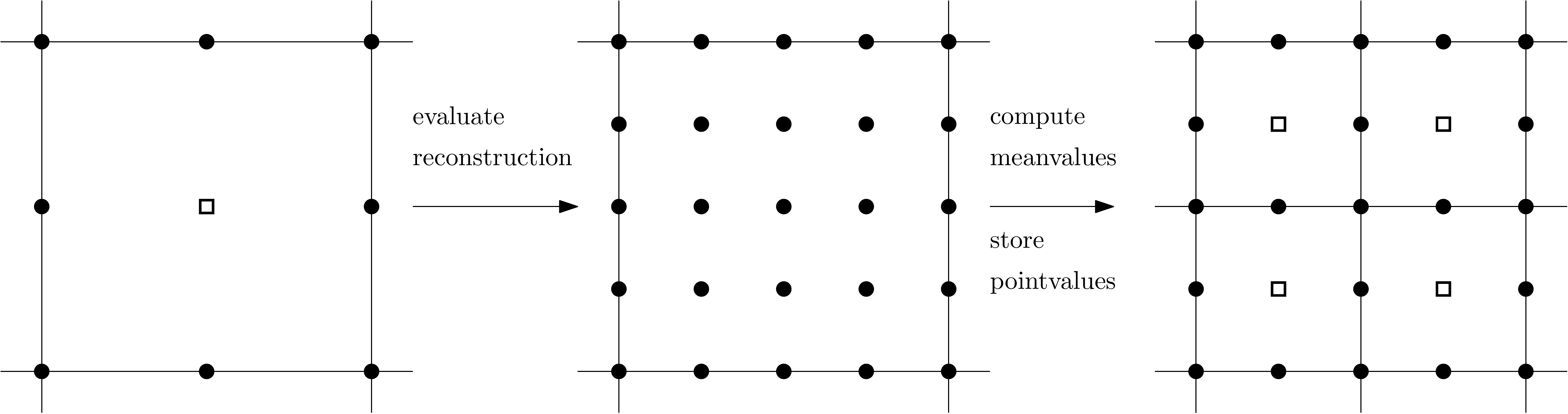}} 
	\caption{\label{fig:transfer2}Illustration of refinement: Point
        values and cell average values on the four fine grid cells are
      computed from the reconstruction of the conserved quantities at
      the coarse grid cell.}
\end{figure}

We use the degrees of freedom of the coarse grid cell to reconstruct a quadratic polynomial as described in Section \ref{sec:AF}. This polynomial can be evaluated at all the required point values along the edges of the fine grid cells. The cell average values of the four fine grid cells are computed using Simpson's rule. This requires the additional computation of four point values on the coarse grid cell at the positions of the  centers of the fine grid cells. An illustration is shown in Figure \ref{fig:transfer2}.   
Since Simpson's rule is exact for our reconstructed function the sum of the cell average values on the four fine grid cells agrees with the cell average on the coarse grid cell exactly. This is a necessary property in order to maintain conservation.

\subsection{Subcycling for a local refinement in space and time}
In explicit finite volume methods for hyperbolic conservation laws, the time step is necessarily restricted by a CFL condition, which requires that the numerical domain of dependence must contain the true domain of dependence of the partial differential equations \cite{book:RJL2002}.

The Active Flux method has a very compact stencil as explained in Section \ref{sec:AF}. Thus, for stability it is necessary to restrict the time step in such a way that information travels at most through one grid cell. Our results from \cite{article:CHK2021} show that the time step should be restricted by $\mbox{CFL} \le 0.75$ for  two-dimensional advection problems and by $\mbox{CFL} \le 0.5$ for the two-dimensional acoustic equations.

On an adaptively refined mesh, the smallest grid cells would typically dictate the time step restriction for the whole domain. To increase the efficiency of the computation, local time stepping (or "subcycling") can be used. In a subcycled computation, several time steps on more refined patches are taken for one time step on the coarsest grid. Subcycling was included in the original AMR algorithm by Berger   and Oliger \cite{be-ol:1984} and Berger and Colella \cite{article:BC1989} and is a standard feature of many AMR codes, including AMRClaw \cite{article:BL1998}, AMReX \cite{zh-al-be-be-bl-ch-da-fr-go:2019} and many others. The local stencil of the Active Flux method allows for efficient implementation of subcycling, which will now be described in more detail. 

\subsubsection*{Subcycling}

In order to use subcycling, we reconstruct the solution in two layers of ghost cells surrounding each patch.  

\begin{figure}[htb]
\centerline{	\includegraphics[width=1.0\textwidth]{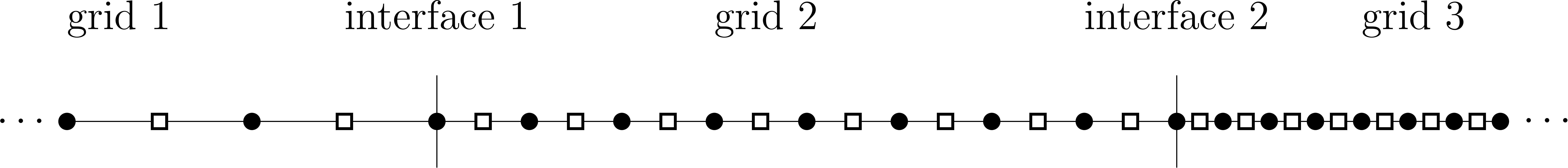} }
\caption{\label{fig:subcycling1}
  Illustration of one-dimensional grid interfaces with different resolution}
\end{figure}

We illustrate the idea for the one-dimensional case but an extension to the two-dimensional situation is straightforward. Assume a situation with three different grid patches as illustrated in Figure \ref{fig:subcycling1}.
The coarsest grid, denoted as grid 1, has refinement level $\ell$. The grid cells on this part of the domain have length $\Delta x_\ell$. Our grid 2 has refinement level $\ell +1$ and the grid cell length satisfies $\Delta x_{\ell +1} = \Delta x_\ell / 2$.    On the finest mesh, i.e. grid 3, we have $\Delta x_{\ell+2} = \Delta x_{\ell +1} / 2$. For efficiency, we will use a time step $\Delta t_\ell$ on grid 1, two time steps  $\Delta t_{\ell+1} = \Delta t_{\ell}/2$ on grid 2 and four time steps $\Delta t_{\ell+2} = \Delta t_{\ell+1}/2 = \Delta t_{\ell}/4$ on grid 3.

In the ForestClaw implementation, we start the time stepping on the finest level grids.  Referring to grids 1,2 and 3 described above, the algorithm proceeds as follows for this three-level mesh configuration.  

\begin{enumerate}
\item Advance the solution one step on grid 3 using time step $\Delta t_{\ell+2}$. 
\item Recursively advance the solution one step on grid 2 (using time step $\Delta t_{\ell+1}$) and on grid 1 (using time step $\Delta t_{\ell}$). 
\item Advance the solution a second step on grid 3.  
\item Grids 2 and 3 are now time synchronized and ghost cell data is exchanged between these levels.
\item Advance the solution a third step on grid 3.  
\item Recursively update the solution a second step on level 2
\item Advance the solution a fourth step on grid 3
\end{enumerate}

\begin{figure}[bht!]
\begin{center}
\includegraphics[width=0.45\textwidth]{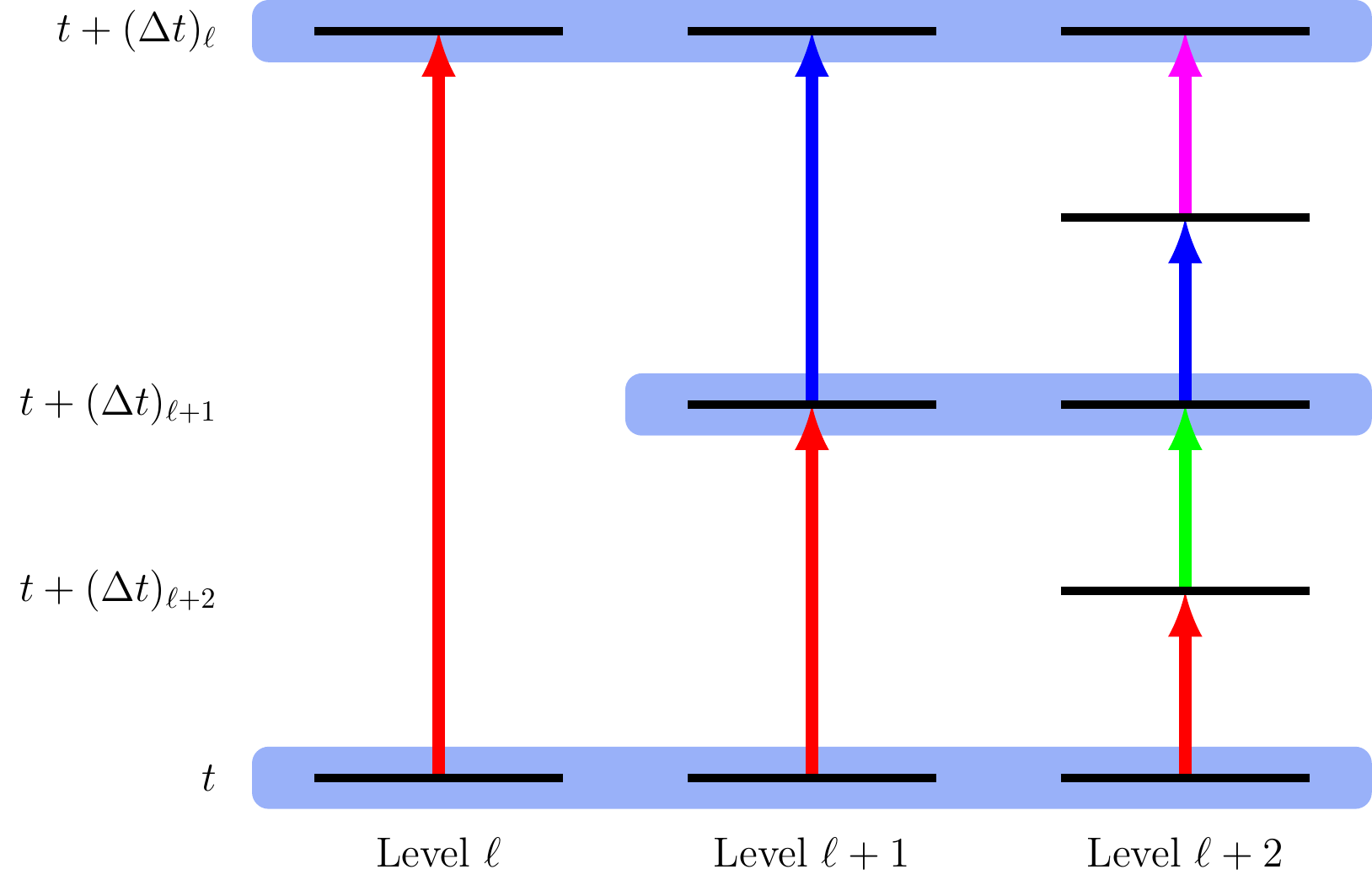}
\end{center}
\caption{Subcycling algorithm for advancing three levels (with grids 1, 2 and 3) a single coarse grid time step of size $\Delta t_{\ell}$.  The algorithm starts at the lower right by advancing one step on the finest level.  All arrows of the same color are advanced recursively. The blue shaded boxes indicate time synchronized levels where ghost exchanges between neighboring grids takes place.  Grids without time synchronized coarse grid neighbors advance by using updated values in an extra layer of ghost cells.}
\label{fig:subcycle}
\end{figure}

In each grid advance, the first layer of ghost cells is updated along with all interior cells.  These ghost cell values are needed to update interior cells at intermediate time levels (e.g. fine time levels that do not exist on coarser levels).   However, whenever two levels are time synchronized, the updated ghost cell data is replaced by data averaged or interpolated from the neighboring finer or coarser grids at the same time level.  For those grids at the physical boundary, physical boundary conditions are used at all time levels.  The time step on each grid uses a stable time step appropriate for that grid.  These steps are illustrated in \Fig{subcycle}.

The key difference between the subcycling approach described above and the subcycling described in the original second order finite volume schemes described by Berger, Oliger and Colella is that in the Active Flux method, ghost cells at the intermediate time level are not filled using time interpolation from neighboring coarser grid cells. Instead we obtain all the required data directly from previously updated degrees of freedom of the Active Flux method.
  
\subsubsection*{Global conservation property}
To ensure the global conservation property of the method, the fluxes at grid cell interfaces need to be defined in a unique way. At regular grid cell interfaces the Active Flux method automatically satisfies this conservation property. At interfaces of patches with different refinement level, conservation needs to be enforced by some kind of ``conservative fix''.  We use the classical approach described by Berger and Colella \cite{article:BC1989}, and update both the coarse and the fine grid cells using the fluxes that have been computed for the more resolved grid.  

\section{Numerical results}
\label{sec:4}
In this section we show numerical results for advection, advective transport with spatially and temporally varying velocity field, Burgers' equation and acoustics. We introduce new Active Flux methods for advective transport problems and discuss the preservation of constant states.

For all computations we use subcycling and the conservative fix described above unless otherwise noted. Each grid patch uses $16 \times 16$ Cartesian grid cells plus ghost cells. Our results confirm third order accuracy of the Active Flux method on adaptively refined grids.

\subsection{Convergence study for advection}
\label{sec:Numerics-advection}
We consider the advection equation (\ref{eqn:advection2d}) on the domain $[0,1]\times[0,1]$ with initial
condition
\begin{equation}\label{eqn:AdvH}
q(x,y,0) = H(r(x,y)+r_0)-H(r(x,y)-r_0),
\end{equation}
where $r(x,y):=\sqrt{(x-x_0)^2+(y-y_0)^2}$ and $H(r):=(\tanh(r/0.02) + 1)/2$ with $r_0=0.15$ and $x_0=y_0=0.5$. We use the advection speeds $a=1$ and $b=0.5$ and  time steps which satisfy $\mbox{CFL}=0.6$.

\begin{figure}[htb]
\begin{center}
	\includegraphics[width=0.35\textwidth]{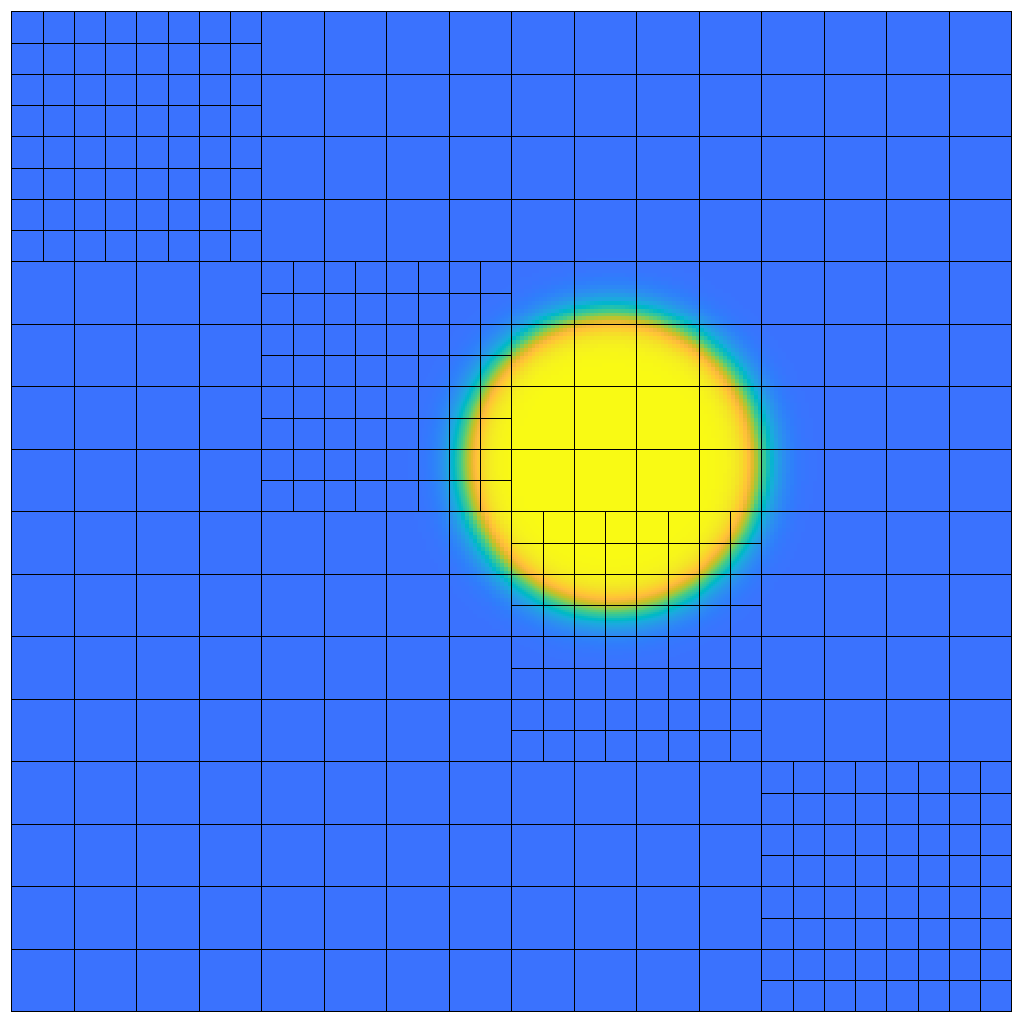}\hfil
	\includegraphics[width=0.35\textwidth]{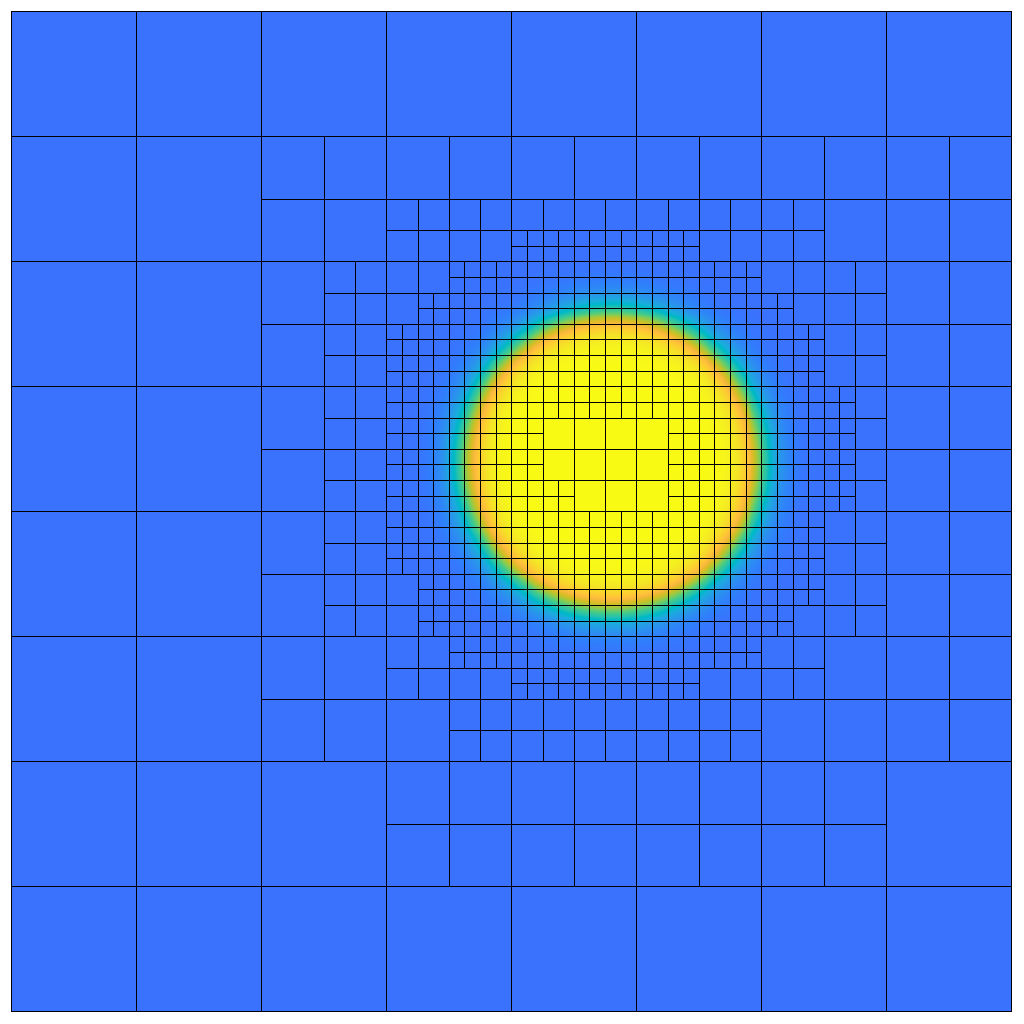}
\end{center}
	\caption{Solutions to (\ref{eqn:AdvH}) with constant refinement along the diagonal (left) and with refinement level 3-6 (right)at time $t=0.1$.} 
  \label{Advection}
\end{figure}

Two different adaptively refined grids are considered.  In the first case refinement is allowed only along the diagonal of the domain as shown in Figure \ref{Advection} (left plot).  In the second case a patch is refined, if $q_{max}-q_{min}>0.001$. In this case the refined grid will follow the solution structure as shown in Figure \ref{Advection} (right plot).  By comparing the numerical solution with the exact solution we can measure the error and compute the experimental order of convergence (EOC).  The results are shown in Table \ref{table:example1} for refinement along the diagonal and \ignore{in Table \ref{table:example1}} for dynamic refinement which follows the solution structure.

\begin{table}[!ht]
  \caption{Error at time $t=1$ measured in the 1-norm
    and convergence rates for the advection problem with constant refinement along the diagonal (left table) and refinement that resolves the relevant solution structure (right table).}
\sisetup{
scientific-notation = true, 
round-mode=places
}
\begin{center}
\begin{minipage}{0.8\textwidth}  
\begin{tabular}{
    *1{S[table-column-width=1.0cm,table-text-alignment=center,round-precision=2]}
    *1{S[table-column-width=1.75cm,table-text-alignment=center,round-precision=2]}
    *1{S[table-column-width=1.5cm,table-text-alignment=center,round-precision=4]}
}
\toprule
{Level} & {Error}         & {EOC}    \\ 
\midrule
{4}     & \num{6.461135e-05}  & {---}    \\ 
{5}     & \num{8.450336e-06}  & \num{2.9347} \\ 
\midrule
{4-5}   & \num{5.075896e-05}  & {---}    \\ 
{5-6}   & \num{6.607801e-06}  & \num{2.9414} \\            
\bottomrule
\end{tabular}\hfill
\begin{tabular}{
    *1{S[table-column-width=1.0cm,table-text-alignment=center,round-precision=2]}
    *1{S[table-column-width=1.75cm,table-text-alignment=center,round-precision=2]}
    *1{S[table-column-width=1.5cm,table-text-alignment=center,round-precision=4]}
}
\toprule
    {Level} & {Error}        & {EOC}    \\ 
\midrule
{5}     & \num{8.450336e-06} & {---} \\ 
{6}     & \num{1.066256e-06} & \num{2.9865} \\
\midrule
{3-5}   & \num{8.467527e-06} & {---} \\ 
{3-6}   & \num{1.069544e-06} & \num{2.9849} \\
\bottomrule
\end{tabular}
\end{minipage}
\end{center}
\label{table:example1}
\end{table}

For the academic test case with refinement along the diagonal, we see that the changes in the grid structure  did not introduce any grid-induced artifacts. The accuracy observed on the adaptively refined grid is comparable with the accuracy on a regular Cartesian grid on the coarser level. 

If the adaptive mesh follows the solution structure, then the accuracy obtained on the adaptively refined grid compares well with the accuracy obtained on a regular Cartesian grid that uses the highest level of refinement in the full domain. These results are shown in Table \ref{table:example1}. The accuracy obtained on the adaptive mesh with levels $3-5$ or $3-6$ compares well with the accuracy obtained on the grids that are refined uniformly to levels $5$ or $6$.  

\subsection{Solid-body rotation}\label{sec:solidBody}
Next we consider the linear advection equation with variable coefficients
\begin{equation}
\partial_t q + \partial_x (a(x,y)q) + \partial_y (b(x,y)q) = 0 \label{eqn:solidbody}
\end{equation}
on $[-1,1] \times [-1,1]$.
We define the advection field for solid body rotation using the 
stream-function $\Psi(x,y)=\frac{\pi}{2}(x^2+y^2)$ and obtain
\begin{equation}
\begin{split}
a(x,y)&=-\partial_y \Psi(x,y) = - \pi y  , \\
b(x,y)&=\phantom{-}\partial_x \Psi(x,y) = \phantom{-} \pi x.
\end{split} \label{eqn:solidbodyODE}
\end{equation}
The velocity field is divergence free, thus Equation
(\ref{eqn:solidbody}) is equivalent to
\begin{equation}
\partial_t q + a(x,y)\partial_x q + b(x,y)\partial_y q = 0. \label{eqn:solidbody2}
\end{equation}
Furthermore, the solution at time $t=2n$, $n \in \mathbb{N}$ agrees with the initial values since the flow has simply made $n$ complete rotations.  The characteristics satisfy the ordinary differential equation
\begin{equation}
\label{eqn:characteristics_adv}
  \begin{aligned}
    x'(t) & = a(x(t),y(t)) \\
    y'(t) & = b(x(t),y(t))
  \end{aligned}
\end{equation}
with initial values $x(0) = x_0$, $y(0) = y_0$. We are interested in $(x(-\tau),y(-\tau))$ with $\tau = \Delta t/2, \Delta t$. Those values can be computed using the analytical solution
\begin{equation}
\label{eqn:characteristics_advtwo}
  \begin{split}
    x(t) & =  - y_0 \sin(\pi t) + x_0 \cos(\pi t) \\
    y(t) & = \phantom{-} y_0  \cos(\pi t) + x_0 \sin(\pi t).
  \end{split}
\end{equation}
This allows us to evaluate the conserved quantities at all required nodes of Simpson's quadrature formula by tracing back the characteristics. For the flux $F_{i+\frac{1}{2},j}$, we use Equation
(\ref{eqn:Simpson-2d}) with, for example,
\begin{equation}
  \begin{aligned}
    f(Q_{i+\frac{1}{2},j}^{n+\frac{1}{2}}) & =
    a(x_{i+\frac{1}{2}},y_j) Q_{i+\frac{1}{2},j}^{n+\frac{1}{2}} =
    -\pi y_j  Q_{i+\frac{1}{2},j}^{n+\frac{1}{2}}\\
    g(Q_{i,j+\frac{1}{2}}^{n+\frac{1}{2}}) & =
    b(x_i,y_{j+\frac{1}{2}}) Q_{i,j+\frac{1}{2}}^{n+\frac{1}{2}} = \pi
    x_i Q_{i,j+\frac{1}{2}}^{n+\frac{1}{2}}
  \end{aligned}
\end{equation}
and
\begin{equation}
  \begin{aligned}
    Q_{i+\frac{1}{2},j}^{n+\frac{1}{2}} & =
    q^n \left( -y_j \sin \left(-\pi \frac{\Delta t}{2}\right) +
    x_{i+\frac{1}{2}}
    \cos\left(-\pi \frac{\Delta t}{2}\right), y_{j} \cos \left(-\pi
    \frac{\Delta t}{2} \right) + x_{i+\frac{1}{2}} \sin \left( -\pi
    \frac{\Delta t}{2} \right) \right) \\
Q_{i,j+\frac{1}{2}}^{n+\frac{1}{2}} & = q^n \left(
  -y_{j+\frac{1}{2}} \sin \left(-\pi \frac{\Delta t}{2}\right) +
    x_{i}
    \cos\left(-\pi \frac{\Delta t}{2}\right), y_{j+\frac{1}{2}} \cos \left(-\pi
    \frac{\Delta t}{2} \right) + x_{i} \sin \left( -\pi
    \frac{\Delta t}{2} \right) \right)
  \end{aligned}
\end{equation}
and analogously for all the other nodes.

We compare numerical solutions of (\ref{eqn:solidbody2}) with initial condition (\ref{eqn:AdvH}) using refinement levels 3-6 and 3-7 after a half rotation with the exact solution and compute the error as well as the EOC. Results are shown in  Table \ref{Tab:solidBodyEOC} (left).
\ignore{
\begin{table}[ht!] 
	\caption{Error at time $t=2$ measured in the 1-norm and EOC for solid-body
          rotation (left) using the approach of Section
          \ref{sec:solidBody} and (right) using the approach of
          Section \ref{sec:swirl}.}
        \hspace*{1cm}
		\begin{tabular}{l|l|l}
			Level & error        & EOC  \\ \hline
                        6   & 4.599185e-06 & \multirow{2}{*}{2.9807} \\ \cline{1-2}
			7   & 5.826258e-07 &   \\ \hline
			3-6   & 4.674628e-06 & \multirow{2}{*}{2.9801} \\ \cline{1-2}
			3-7   & 5.924630e-07 &                        
		\end{tabular} \hfill
                \begin{tabular}{l|l|l}
		Level & error        & EOC                     \\ \hline
		6   & 5.791702e-06 & \multirow{2}{*}{2.9786 } \\ \cline{1-2}
		7   & 7.347671e-07  &   \\ \hline
		3-6   & 5.792193e-06  & \multirow{2}{*}{2.9781} \\ \cline{1-2}
		3-7   & 7.351174e-07 &                        
	\end{tabular}    \hspace*{1cm}          
	 \label{Tab:solidBodyEOC}
\end{table}
}

\begin{table}[ht!] 
  \caption{Error at time $t=2$ measured in the 1-norm and EOC for solid-body
          rotation (left) using the approach of Section
          \ref{sec:solidBody} and (right) using the approach of
          Section \ref{sec:swirl}.}
\sisetup{
scientific-notation = true, 
round-mode=places, round-precision=2}
\begin{center}
\begin{minipage}{0.8\textwidth}
\begin{tabular}{
    *1{S[table-column-width=1.0cm,table-text-alignment=center,round-precision=2]}
    *1{S[table-column-width=1.75cm,table-text-alignment=center,round-precision=2]}
    *1{S[table-column-width=1.5cm,table-text-alignment=center,round-precision=4]}
}
\toprule
      {Level} & {Error}        & {EOC}  \\ 
\midrule
      {6}   & \num{4.599185e-06} & {---} \\ 
      {7}   & \num{5.826258e-07} &   \num{2.9807}\\ 
\midrule      
      {3-6}   & \num{4.674628e-06} & {---} \\ 
      {3-7}   & \num{5.924630e-07} & \num{2.9801} \\
\bottomrule              
    \end{tabular} 
    \hfill
\begin{tabular}{
    *1{S[table-column-width=1.0cm,table-text-alignment=center,round-precision=2]}
    *1{S[table-column-width=1.75cm,table-text-alignment=center,round-precision=2]}
    *1{S[table-column-width=1.5cm,table-text-alignment=center,round-precision=4]}
}
\toprule 
{Level} & {Error}        & {EOC}  \\ 
\midrule    
    {6}   & \num{5.791702e-06} & {---} \\ 
    {7}   & \num{7.347671e-07}  &   \num{2.9786} \\
\midrule    
    {3-6}   & \num{5.792193e-06}  & {---} \\ 
    {3-7}   & \num{7.351174e-07} &  \num{2.9781} \\ 
\bottomrule                               
  \end{tabular}
\end{minipage}
\end{center}  
   \label{Tab:solidBodyEOC}
\end{table}

\begin{figure}[ht!]
\begin{center}
\includegraphics[width=0.35\textwidth]{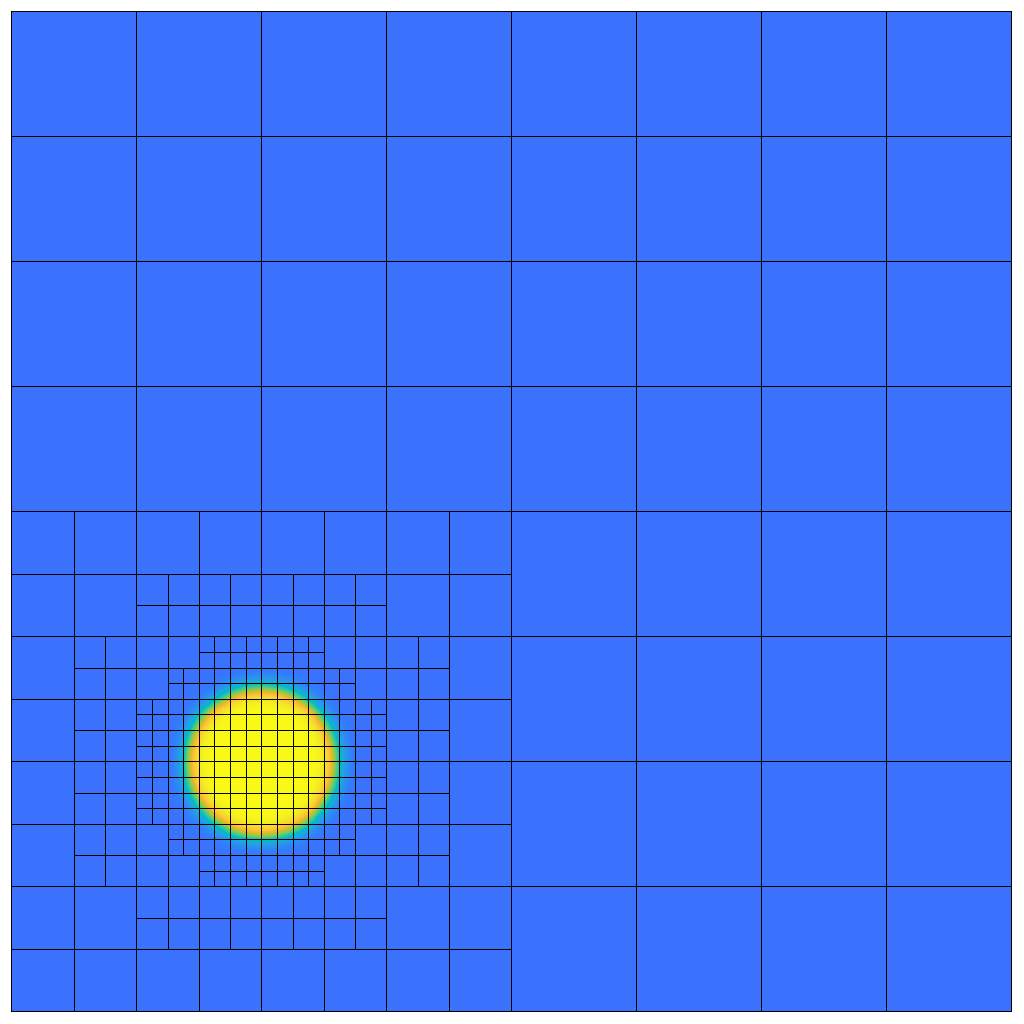}\hfil
\includegraphics[width=0.35\textwidth]{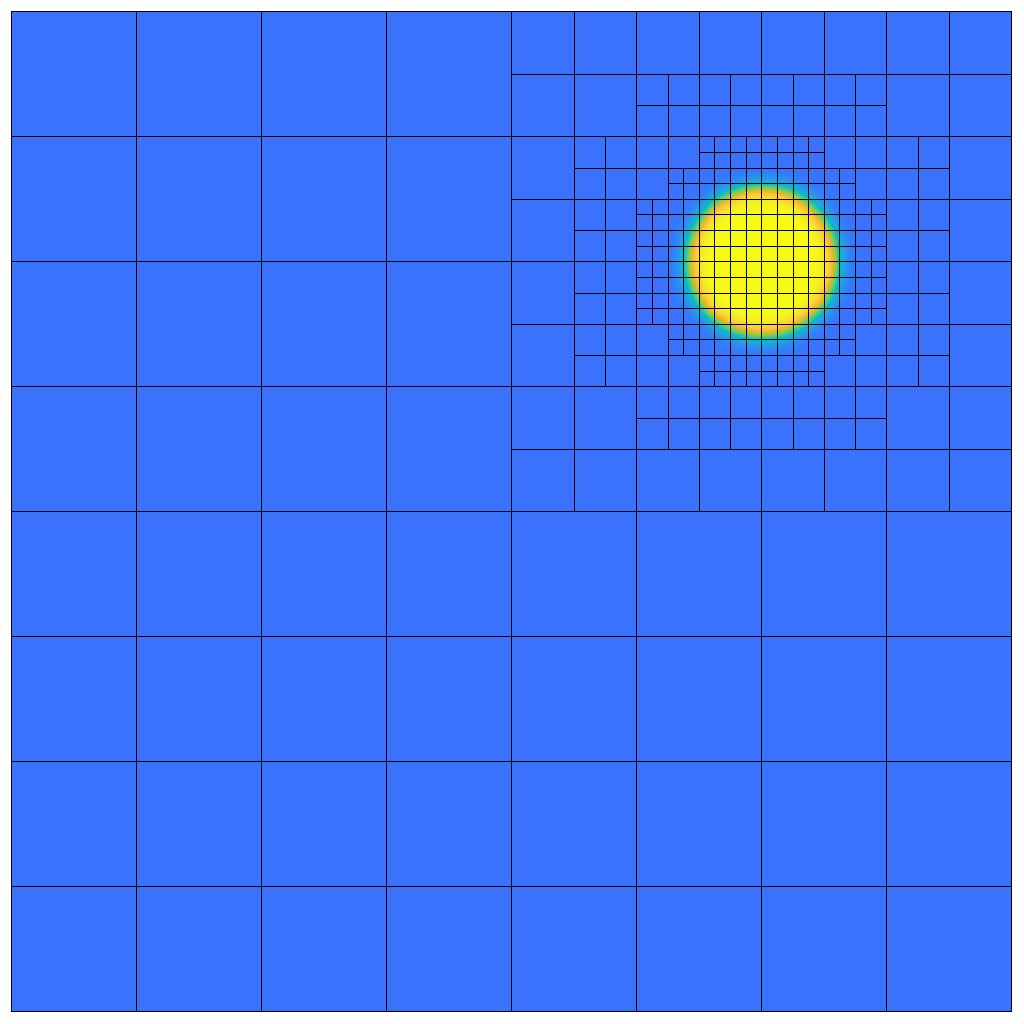}
\caption{Solution to (\ref{eqn:solidbody2}) with initial
          condition (\ref{eqn:AdvH}) after a half (left) and a full
          rotation (right).}
\label{solidBody}
\end{center}
\end{figure}

Again, the accuracy of the computations on the  adaptively refined grids compares well with the accuracy obtained on the equidistant grids with highest resolution.  Figure \ref{solidBody} shows the numerical solution with refinement level 3-6 after a half  and a full rotation. 

\begin{remark}
The Active Flux method for solid body rotation described in Section
\ref{sec:solidBody} preserves constant states on regular Cartesian grids.
\end{remark}
\textbf{Proof:} We consider constant data at time $t_n$, i.e.\
$q(x,y,t_n) = C \in \mathbb{R}$. Then all the point values at time $t_n$,
$t_{n+\frac{1}{2}}$ and $t_{n+1}$ are equal to $C$ and the numerical fluxes have
the form
\begin{equation*}
  \begin{split}
F_{i+\frac{1}{2},j} & = - \frac{C \pi}{6} \left(y_{j-\frac{1}{2}} + 4 y_j +
  y_{j+\frac{1}{2}} \right) \\
G_{i,j+\frac{1}{2}} & = \frac{C \pi}{6} \left( x_{i-\frac{1}{2}} + 4 x_i
  + x_{i+\frac{1}{2}} \right).
\end{split}
\end{equation*}
Thus, the fluxes in the $x$-direction only depend on $y$ and the fluxes in the $y$-direction only depend on $x$. Consequently, all the flux differences cancel and the cell average values remain constant.

\hfill $\square$\\

\subsection{The swirl flow problem}\label{sec:swirl}
Now we consider advective transport in a velocity field that depends on space and time.  For this case, we use
\begin{equation}\label{eqn:advection_spacetime}
\partial_t q + \partial_x \left(a(x,y,t) q\right) + \partial_y \left(
  b(x,y,t) q \right) = 0.
\end{equation}
The velocity field is again defined via a stream function, which now has the form
$$
\Psi(x,y,t) = \frac{1}{\pi} \sin^2(\pi x) \cdot \sin^2(\pi y)\cdot
\cos(\pi t).
$$
We use finite difference approximations to compute
\begin{equation}
  \begin{split}
a(x,y,t)&=-\frac{\Psi(x,y+\Delta y/2,t)-\Psi(x,y-\Delta y/2,t)}{\Delta y}, \\
b(x,y,t)&= \phantom{-}\frac{\Psi(x+\Delta x/2,y,t)-\Psi(x-\Delta x/2,y,t)}{\Delta x}.
\end{split} \label{eqn:swirlODE}
\end{equation}
We trace the characteristics back in space and time using the classical fourth order accurate Runge-Kutta method backwards in time.

For the flux computation of $F_{i+\frac{1}{2},j}$ we use the approximation
\begin{equation}\label{eqn:flux_swirl}
  \begin{split}
    & \frac{1}{\Delta t \Delta y}
    \int_{t_n}^{t_{n+1}} \int_{y_{j-\frac{1}{2}}}^{y_{j+\frac{1}{2}}}
    a(x_{i+\frac{1}{2}},y,t) q(x_{i+\frac{1}{2}},y,t) dy dt \\
    & \approx  \frac{1}{\Delta t \Delta y}
    \int_{t_n}^{t_{n+1}} \left( a(x_{i+\frac{1}{2}},y_j,t)
  \int_{y_{j-\frac{1}{2}}}^{y_{j+\frac{1}{2}}} q(x_{i+\frac{1}{2}},y,t)
  dy \right) dt
\end{split}
\end{equation}
and analogously for the flux $G_{i,j+\frac{1}{2}}$. Using again Simpson's
rule we obtain
\begin{equation}\label{eqn:swirlF}
  \begin{split}
F_{i+\frac{1}{2},j} & = \frac{1}{6} \Bigg{(}
 \frac{1}{6} a(x_{i+\frac{1}{2}},y_j,t_n)  \cdot \left(
   Q_{i+\frac{1}{2},j-\frac{1}{2}}^n + 4 Q_{i+\frac{1}{2},j}^n +
   Q_{i+\frac{1}{2},j+\frac{1}{2}}^n \right) \\
 & \hspace*{1cm} + \frac{4}{6} a(x_{i+\frac{1}{2}},y_j,t_{n+\frac{1}{2}}) \cdot
 \left( Q_{i+\frac{1}{2},j-\frac{1}{2}}^{n+\frac{1}{2}} + 4
   Q_{i+\frac{1}{2},j}^{n+\frac{1}{2}} +
   Q_{i+\frac{1}{2},j+\frac{1}{2}}^{n+\frac{1}{2}} \right) \\
  & \hspace*{1cm} + \frac{1}{6} a(x_{i+\frac{1}{2}},y_j,t_{n+1}) \cdot \left(
   Q_{i+\frac{1}{2},j-\frac{1}{2}}^{n+1} + 4 Q_{i+\frac{1}{2},j}^{n+1}
   + Q_{i+\frac{1}{2},j+\frac{1}{2}}^{n+1} \right)  \Bigg{)}
  \end{split}
\end{equation}
and
\begin{equation}\label{eqn:swirlG}
  \begin{split}
G_{i,j+\frac{1}{2}} & = \frac{1}{6} \Bigg{(}
 \frac{1}{6} b(x_{i},y_{j+\frac{1}{2}},t_n)  \cdot \left(
   Q_{i-\frac{1}{2},j+\frac{1}{2}}^n + 4 Q_{i,j+\frac{1}{2}}^n +
   Q_{i+\frac{1}{2},j+\frac{1}{2}}^n \right) \\
 & \hspace*{1cm} + \frac{4}{6} b(x_{i},y_{j+\frac{1}{2}},t_{n+\frac{1}{2}}) \cdot
 \left( Q_{i-\frac{1}{2},j+\frac{1}{2}}^{n+\frac{1}{2}} + 4
   Q_{i,j+\frac{1}{2}}^{n+\frac{1}{2}} +
   Q_{i+\frac{1}{2},j+\frac{1}{2}}^{n+\frac{1}{2}} \right) \\
  & \hspace*{1cm} + \frac{1}{6} b(x_{i},y_{j+\frac{1}{2}},t_{n+1}) \cdot \left(
   Q_{i-\frac{1}{2},j+\frac{1}{2}}^{n+1} + 4 Q_{i,j+\frac{1}{2}}^{n+1}
   + Q_{i+\frac{1}{2},j+\frac{1}{2}}^{n+1} \right)  \Bigg{)}.
  \end{split}
\end{equation}

\begin{theorem}\label{theo:constant}
  The Active Flux method with fluxes of the form (\ref{eqn:swirlF}),
  (\ref{eqn:swirlG}) and  $a, b$ as defined in
  (\ref{eqn:swirlODE}), provides a method for
  (\ref{eqn:advection_spacetime}) that preserves constant states on
  regular Cartesian grids.
\end{theorem}
\textbf{Proof:} We consider constant data at time $t_n$, i.e.\
$q(x,y,t_n)=C \in \mathbb{R}$.
Then all the point values at time $t_{n}$, $t_{n+\frac{1}{2}}$ and
$t_{n+1}$, which are obtained by tracing back the characteristics, are
also equal to $C$. The finite volume update now reduces to 
\begin{equation*}
  \begin{split}
Q_{i,j}^{n+1} = Q_{i,j}^n & - \frac{\Delta t}{\Delta x} \frac{C}{6}
\Big{(}  a(x_{i+\frac{1}{2}},y_j,t_n) + 4
a(x_{i+\frac{1}{2}},y_j,t_{n+\frac{1}{2}}) +
a(x_{i+\frac{1}{2}},y_j,t_{n+1}) \\
& \hspace*{1cm} -a(x_{i-\frac{1}{2}},y_j,t_n) -
4a(x_{i-\frac{1}{2}},y_j,t_{n+\frac{1}{2}}) -
a(x_{i-\frac{1}{2}},y_j,t_{n+1}) \Big{)} \\
& - \frac{\Delta t}{\Delta y} \frac{C}{6} \Big{(}
b(x_i,y_{j+\frac{1}{2}},t_n) + 4
b(x_i,y_{j+\frac{1}{2}},t_{n+\frac{1}{2}})+b(x_i,y_{j+\frac{1}{2}}
t_{n+1}) \\
& \hspace*{1cm} - b(x_i,y_{j-\frac{1}{2}},t_n) - 4
b(x_i,y_{j-\frac{1}{2}},t_{n+\frac{1}{2}}) -b(x_i,y_{j-\frac{1}{2}},t_{n+1}) \Big{)}.
  \end{split}
\end{equation*}
For the terms at time level $t_n$ we obtain
\begin{equation*}
  \begin{split}
& - \frac{C \Delta t}{6 \Delta x}  \left(
a(x_{i+\frac{1}{2}},y_j,t_n) - a(x_{i-\frac{1}{2}},y_j,t_n) \right) -
\frac{C \Delta t}{6 \Delta y} \left( b(x_i,y_{j+\frac{1}{2}},t_n) -
  b(x_i,y_{j-\frac{1}{2}},t_n) \right) \\
& = - \frac{C \Delta t}{6 \Delta x \Delta y} \Big{(} -
  \Psi(x_{i+\frac{1}{2}},y_{j+\frac{1}{2}},t_n) +
  \Psi(x_{i+\frac{1}{2}},y_{j-\frac{1}{2}},t_n) 
 + \Psi(x_{i-\frac{1}{2}},y_{j+\frac{1}{2}},t_n) -
  \Psi(x_{i-\frac{1}{2}},y_{j-\frac{1}{2}},t_n) \Big{)} \\
  & \quad - \frac{C \Delta t}{6  \Delta y \Delta x} \Big{(}
  \Psi(x_{i+\frac{1}{2}},y_{j+\frac{1}{2}},t_n) -
  \Psi(x_{i-\frac{1}{2}},y_{j+\frac{1}{2}},t_n) 
 - \Psi(x_{i+\frac{1}{2}},y_{j-\frac{1}{2}},t_n) +
  \Psi(x_{i-\frac{1}{2}},y_{j+\frac{1}{2}},t_n) \Big{)} \\
  & = 0
  \end{split}
\end{equation*}
In the same way the terms at time $t_{n+\frac{1}{2}}$ and $t_{n+1}$ cancel and we obtain $Q_{i,j}^{n+1} = Q_{i,j}^n = C$ for all $i,j$.

\hfill $\square$\\

\begin{figure}[htb]
	\begin{center}
\includegraphics[width=0.5\textwidth]{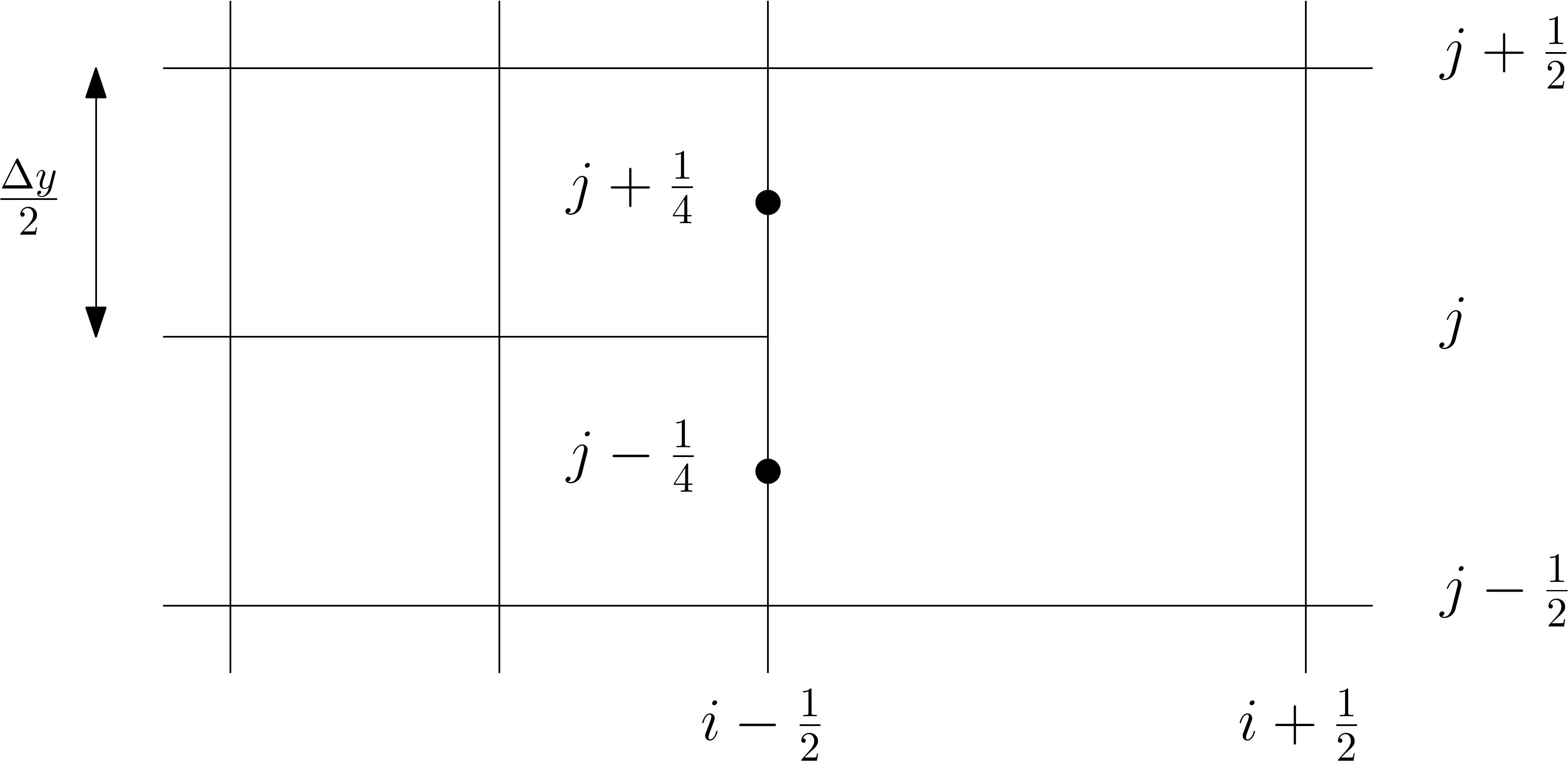}
\end{center}
	\caption{\label{fig:constantStatesAMR} Coarse grid cell with
          neighboring fine grid cells.
}
\end{figure}

Now we consider Cartesian grids with adaptive mesh refinement. Without loss of generality we restrict considerations to the situation shown in Figure \ref{fig:constantStatesAMR}.
We again consider constant data at time $t_n$, i.e.\ $q(x,y,t_n) = C$. The update of the small grid cells preserves by Theorem \ref{theo:constant} constant states. For the large cell update, the left flux $F_{i-\frac{1}{2},j}$ is the sum of two fluxes used to update the two small cells, i.e.\
\begin{equation*}
  \begin{split}
F_{i-\frac{1}{2},j} & = \frac{1}{2} \left(
  F_{i-\frac{1}{2},j-\frac{1}{4}} + F_{i-\frac{1}{2},j+\frac{1}{4}}
\right) \\
& = \frac{1}{2} C \Big{[} \frac{1}{6} \left(
  a(x_{i-\frac{1}{2}},y_{j-\frac{1}{4}},t_n) + 4
  a(x_{i-\frac{1}{2}},y_{j-\frac{1}{4}},t_{n+\frac{1}{2}}) +
  a(x_{i-\frac{1}{2}},y_{j-\frac{1}{4}},t_{n+1}) \right) \\
& \hspace*{1.cm}  + \frac{1}{6} \left(
  a(x_{i-\frac{1}{2}},y_{j+\frac{1}{4}},t_n) + 4
  a(x_{i-\frac{1}{2}},y_{j+\frac{1}{4}},t_{n+\frac{1}{2}}) +
  a(x_{i-\frac{1}{2}},y_{j+\frac{1}{4}},t_{n+1}) \right) \Big{]} \\
& = \frac{1}{2} C \frac{1}{6} \Big{[}
-
\frac{\Psi(x_{i-\frac{1}{2}},y_j,t_n)-\Psi(x_{i-\frac{1}{2}},y_{j-\frac{1}{2}},t_n)}{\Delta
y / 2} - 4 \frac{\Psi(x_{i-\frac{1}{2}},y_j,t_{n+\frac{1}{2}})-\Psi(x_{i-\frac{1}{2}},y_{j-\frac{1}{2}},t_{n+\frac{1}{2}})}{\Delta
y / 2}\\
& - \frac{\Psi(x_{i-\frac{1}{2}},y_j,t_{n+1})-\Psi(x_{i-\frac{1}{2}},y_{j-\frac{1}{2}},t_{n+1})}{\Delta
y / 2} - \frac{\Psi(x_{i-\frac{1}{2}},y_{j+\frac{1}{2}},t_{n})-\Psi(x_{i-\frac{1}{2}},y_{j},t_{n})}{\Delta
y / 2} \\
& -4 \frac{\Psi(x_{i-\frac{1}{2}},y_{j+\frac{1}{2}},t_{n+\frac{1}{2}})-\Psi(x_{i-\frac{1}{2}},y_{j},t_{n+\frac{1}{2}})}{\Delta
y / 2} -
\frac{\Psi(x_{i-\frac{1}{2}},y_{j+\frac{1}{2}},t_{n+1})-\Psi(x_{i-\frac{1}{2}},y_{j},t_{n+1})}{\Delta
  y / 2} \Big{]}\\
& = \frac{1}{6} C \Big{[}
- \frac{\Psi(x_{i-\frac{1}{2}},y_{j+\frac{1}{2}},t_n) -
  \Psi(x_{i-\frac{1}{2}},y_{j-\frac{1}{2}},t_n)}{\Delta y}  \\
& - 4 \frac{\Psi(x_{i-\frac{1}{2}},y_{j+\frac{1}{2}},t_{n+\frac{1}{2}}) -
  \Psi(x_{i-\frac{1}{2}},y_{j-\frac{1}{2}},t_{n+\frac{1}{2}})}{\Delta
  y} - \frac{\Psi(x_{i-\frac{1}{2}},y_{j+\frac{1}{2}},t_{n+1}) -
  \Psi(x_{i-\frac{1}{2}},y_{j-\frac{1}{2}},t_{n+1})}{\Delta
  y} \Big{]}\\
& = C \frac{1}{6} \left( a(x_{i-\frac{1}{2}},y_j,t_n) + 4
  a(x_{i-\frac{1}{2}},y_j,t_{n+\frac{1}{2}}) +
  a(x_{i-\frac{1}{2}},y_j,t_{n+1}) \right).
  \end{split}
\end{equation*}
This is the same flux that would be used on a regular Cartesian grid. Thus, the update of the small grid cells as well as the update of the coarse grid cells preserve by Theorem \ref{theo:constant} constant states. We summarize our result.
\begin{corollary}
The Active Flux method described above preserves constant states on Cartesian grids with adaptive mesh refinement \ignore{as long as no subcycling is used.} when used with global time stepping (e.g. no subcycling). 
\end{corollary}
With subcycling, the sum of the fluxes from the fine grid cells and that is used to update the coarse grid cell would contain components at the intermediate time that are in general not balanced by the remaining coarse grid fluxes.

To perform a numerical convergence study, we consider smooth initial values of the form
\begin{equation}\label{eqn:id_swirl}
q(x,y,0)=\exp(-100((x-0.5)^2+(y-0.25)^2))+\exp(-100((x-0.5)^2+(y-0.75)^2))
\end{equation}
on the domain $[0,1]\times[0,1]$ with periodic boundary conditions.
\begin{figure}[htb]
\begin{center}
	\includegraphics[width=0.28\textwidth]{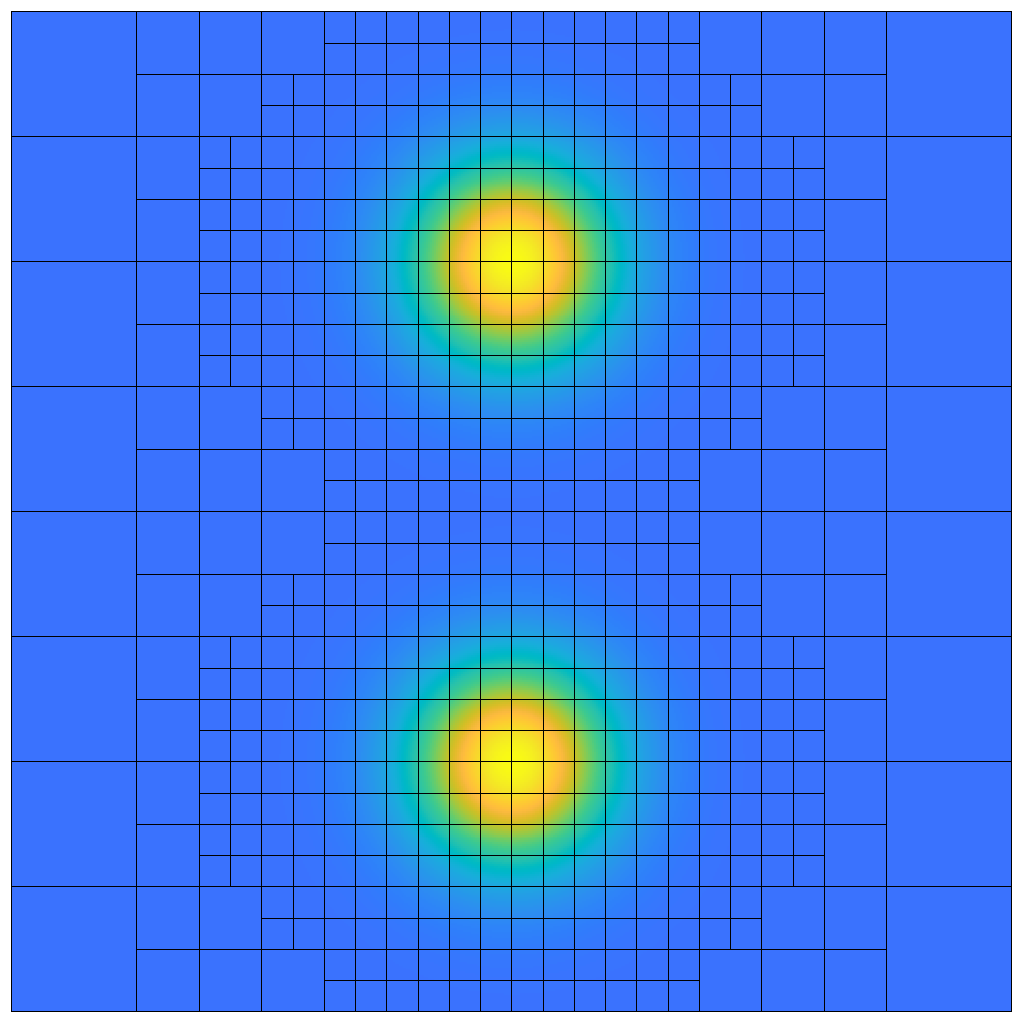}\hfil
	\includegraphics[width=0.28\textwidth]{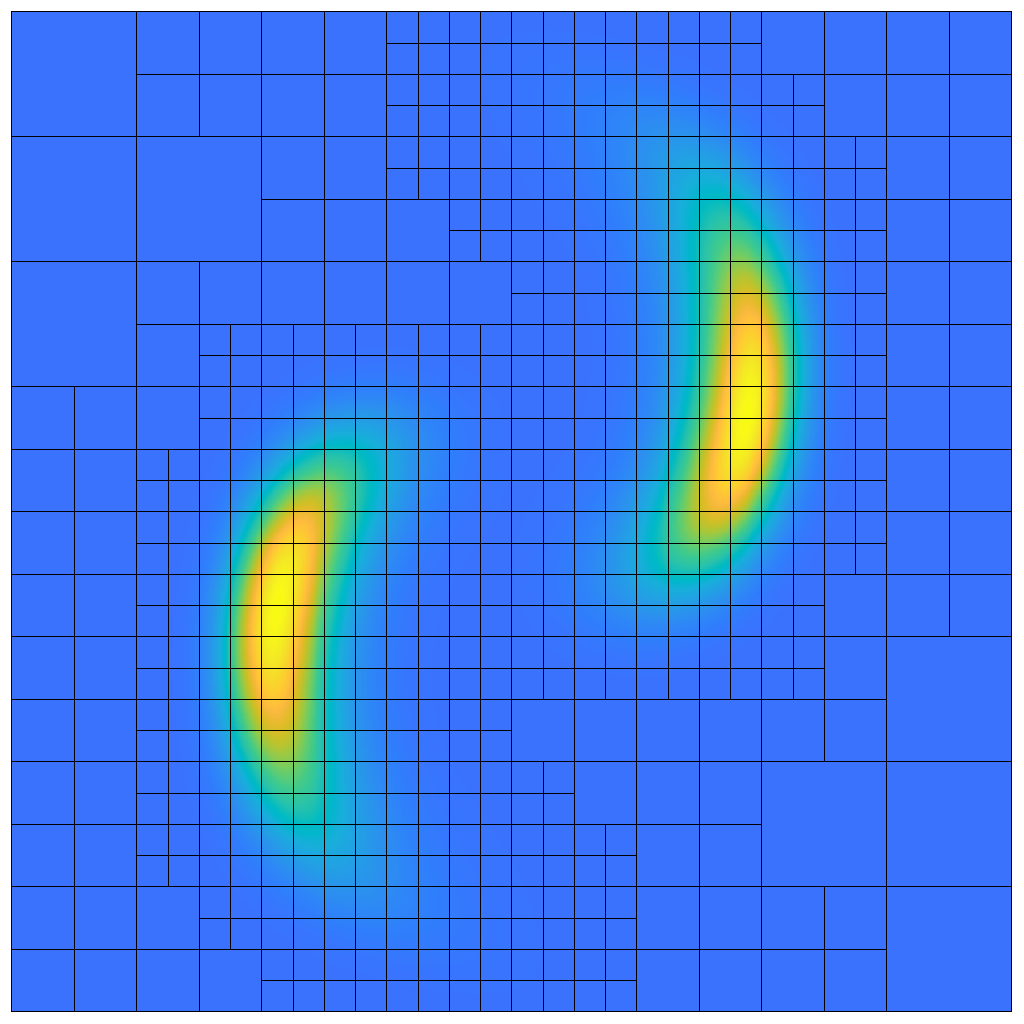}\hfil
	\includegraphics[width=0.28\textwidth]{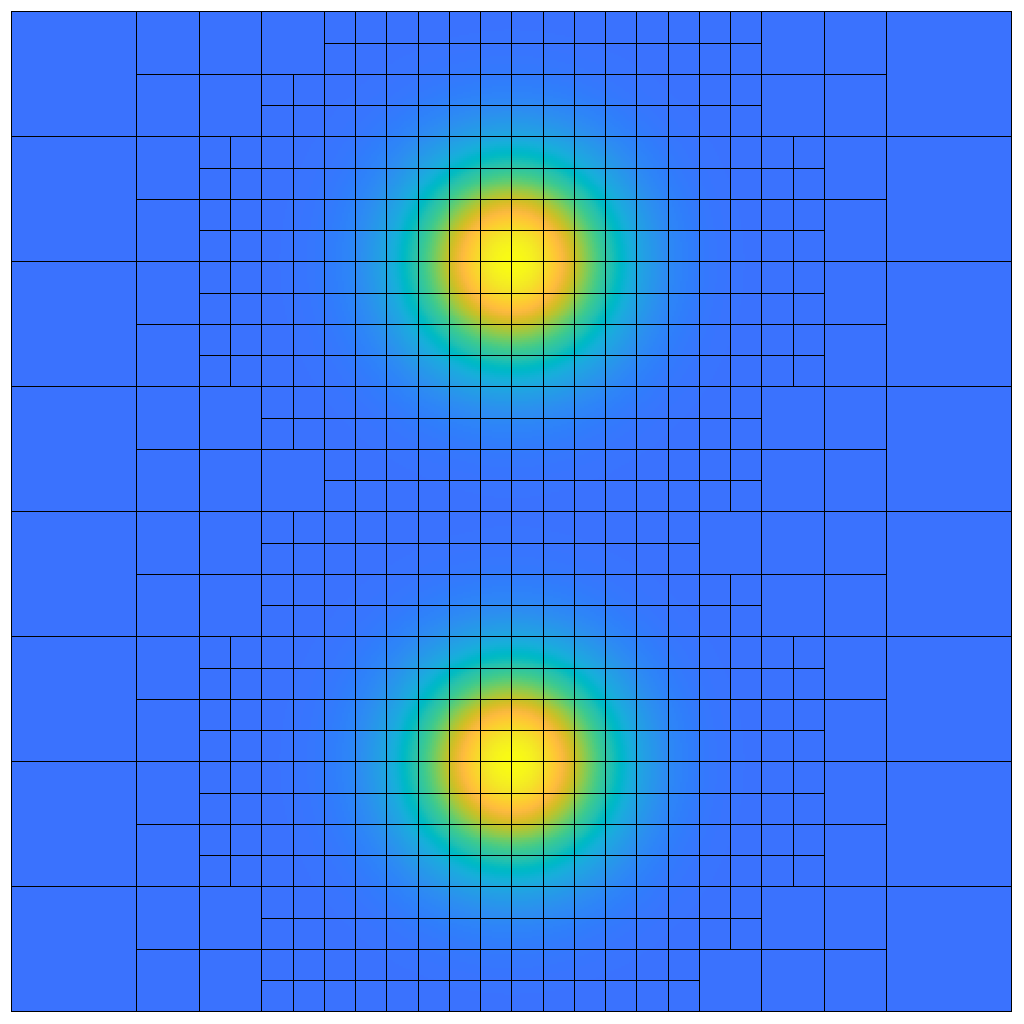}
\end{center}  
	\caption{Solution of the swirl problem with initial data of
          the form (\ref{eqn:id_swirl}) at times
		$t = 0$ (left),
		$t = 0.5$ (middle) and 
		$t = 1$ (right) for the swirl problem with smooth initial condition.} \label{swirlSmooth}
\end{figure}
Numerical results are shown in Figure \ref{swirlSmooth} and the results of  a convergence study with and without subcycling is shown in Table \ref{Tab:swirlEOC}. The results again confirm third order convergence.  The use of subcycling did not significantly influence the accuracy.
\begin{table}[htb!] 
  \caption{Error at time $t=1$ measured in the 1-norm and EOC for the smooth swirl
    problem  (left) with subcycling, (right) without subcycling.}
\sisetup{
scientific-notation = true, 
round-mode=places, round-precision=2}
\begin{center}
  \begin{minipage}{0.8\textwidth}
\begin{tabular}{
    *1{S[table-column-width=1.0cm,table-text-alignment=center,round-precision=2]}
    *1{S[table-column-width=1.75cm,table-text-alignment=center,round-precision=2]}
    *1{S[table-column-width=1.5cm,table-text-alignment=center,round-precision=4]}
}
\toprule
      {Level} & {Error}        & {EOC}  \\ 
  \midrule
{4}     & \num{5.722927e-06} & {---} \\ 
  {5}   & \num{7.268088e-07} &  \num{2.9771}  \\ 
\midrule
{3-4}   & \num{5.758274e-06} & {---} \\ 
  {3-5}   & \num{7.346642e-07} & \num{2.9705}  \\                   
\bottomrule
\end{tabular}\hfill
\begin{tabular}{
    *1{S[table-column-width=1.0cm,table-text-alignment=center,round-precision=2]}
    *1{S[table-column-width=1.75cm,table-text-alignment=center,round-precision=2]}
    *1{S[table-column-width=1.5cm,table-text-alignment=center,round-precision=4]}
}
\toprule
      {Level} & {Error}        & {EOC}  \\ 
  \midrule
{4}     & \num{5.722927e-06} & {---} \\ 
  {5}   & \num{7.268088e-07} &  \num{2.9771}  \\ 
\midrule
{3-4}   & \num{5.724930e-06} & {---} \\ 
  {3-5}   & \num{7.305236e-07} & \num{2.9703}  \\                   
\bottomrule
\end{tabular}
\end{minipage}
\end{center}
\label{Tab:swirlEOC}
\end{table}

To measure the error in the approximation of constant states by using subcycling, we compute the solution at time $t=1$ using the same velocity field but with constant initial values that are equal to one in the whole domain. Refinement with levels $3-4$ is used along the diagonal as shown in Figure \ref{Advection} (left). With subcycling we observe an error of size $10^{-13}$. Without subcycling the error is about $10^{-16}$, i.e.\ agrees with the expected truncation error. 

\begin{figure}[htb]
\begin{center}
\includegraphics[width=0.28\textwidth]{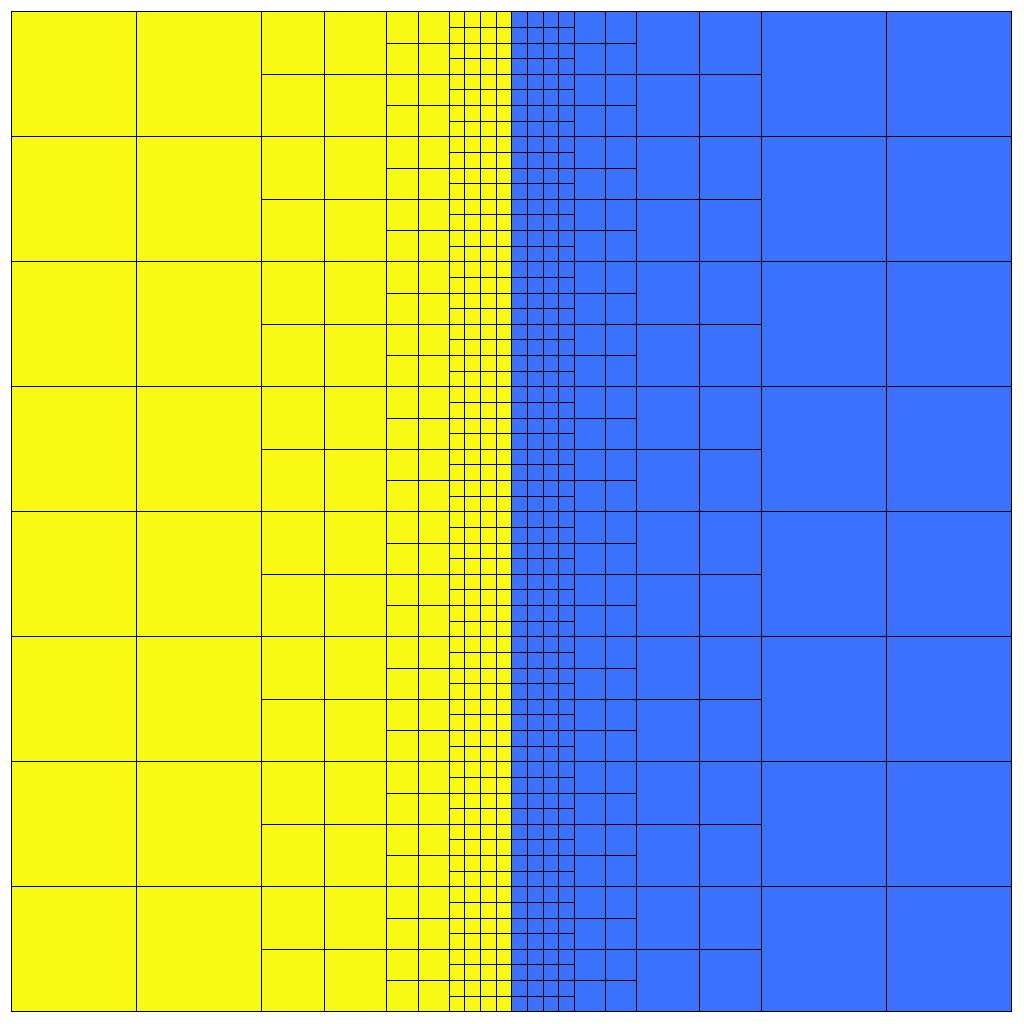}\hfil
\includegraphics[width=0.28\textwidth]{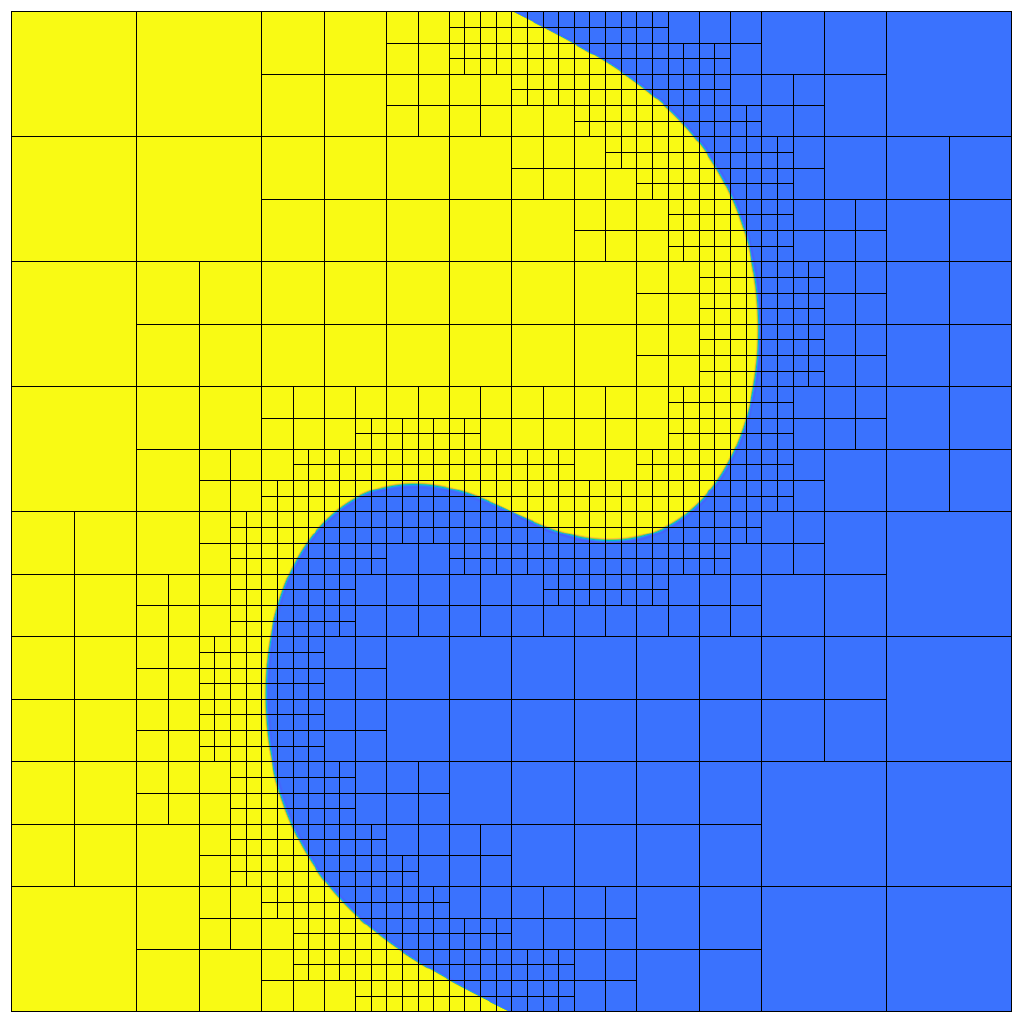}\hfil
\includegraphics[width=0.28\textwidth]{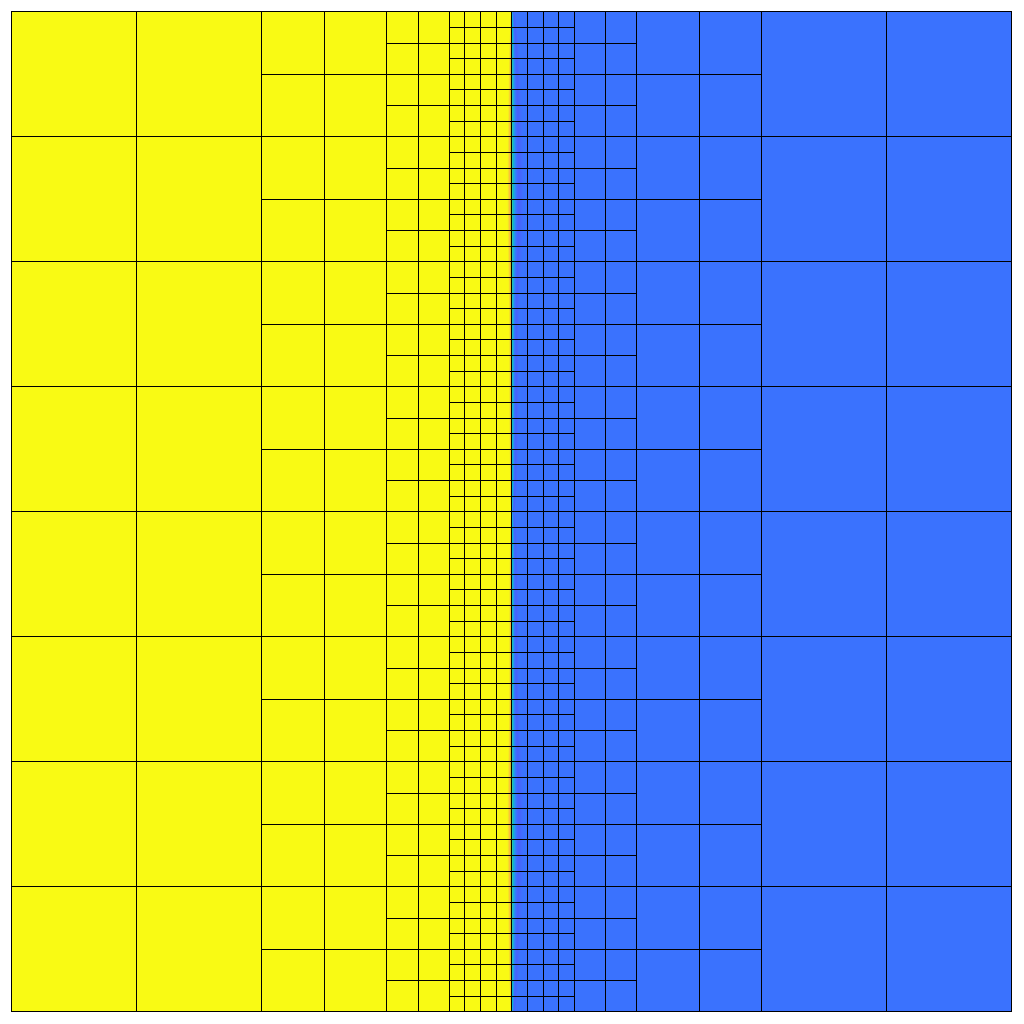}
\end{center}
	\caption{\label{fig:swirlRP} Solution of the swirl problem
          with piecewise constant initial values at times
		$t = 0$ (left),
		$t = 0.5$ (middle) and 
		$t = 1$ (right).} \label{swirlRiemann}
\end{figure}

We now also consider the swirl problem with piecewise constant initial values of the form
\begin{equation*}
q(x,y,0)=\begin{cases} 
	1 & 0 \leq x \leq 0.5 \\
	0 & 0.5 < x \leq 1 
\end{cases} \end{equation*}
on the domain $[0,1]\times[0,1]$ with periodic boundaries on top and bottom and with zero-order extrapolation on the left and right boundary using refinement level 3-6. Results are shown in Figure \ref{fig:swirlRP}.    
Subcycling was used although this slightly perturbs the exact preservation of constant states.
 
While the method from Section \ref{sec:solidBody} preserves constant states for solid body rotation due to the special form of the velocity field, we can also compute fluxes for solid body rotation using the more general constant state preserving approach presented in this section. Numerical results are shown in Table \ref{Tab:solidBodyEOC} (right). The error is slightly larger, due to the use of the approximation described in Equation (\ref{eqn:flux_swirl}), but the method clearly preserves third order accuracy.

\subsection{Burgers' equation}
It is straight forward to apply the adaptive Active Flux method to scalar nonlinear hyperbolic problems. We consider the Burgers' equation (\ref{eqn:burgers2d}) on the domain $[0,1] \times [0,1]$ with initial values
\begin{equation}\label{eqn:sinsin01}
q(x,y,0) = \sin(2\pi x)\sin(2\pi y)+0.1.
\end{equation}

To check the accuracy we compute numerical solutions at an early time, here we use $t=0.05$,  at which the solution structure is still smooth. The time steps satisfy $\mbox{CFL} \le 0.5$. Furthermore, we use subcycling and the conservative fix. To test the accuracy of the adaptive method we enforce refinement along the diagonal from the upper left patch to the lower right patch. This test (not shown here) confirms third order accuracy also for this nonlinear problem.

\begin{figure}[htb]
\begin{center}
\includegraphics[width=0.28\textwidth]{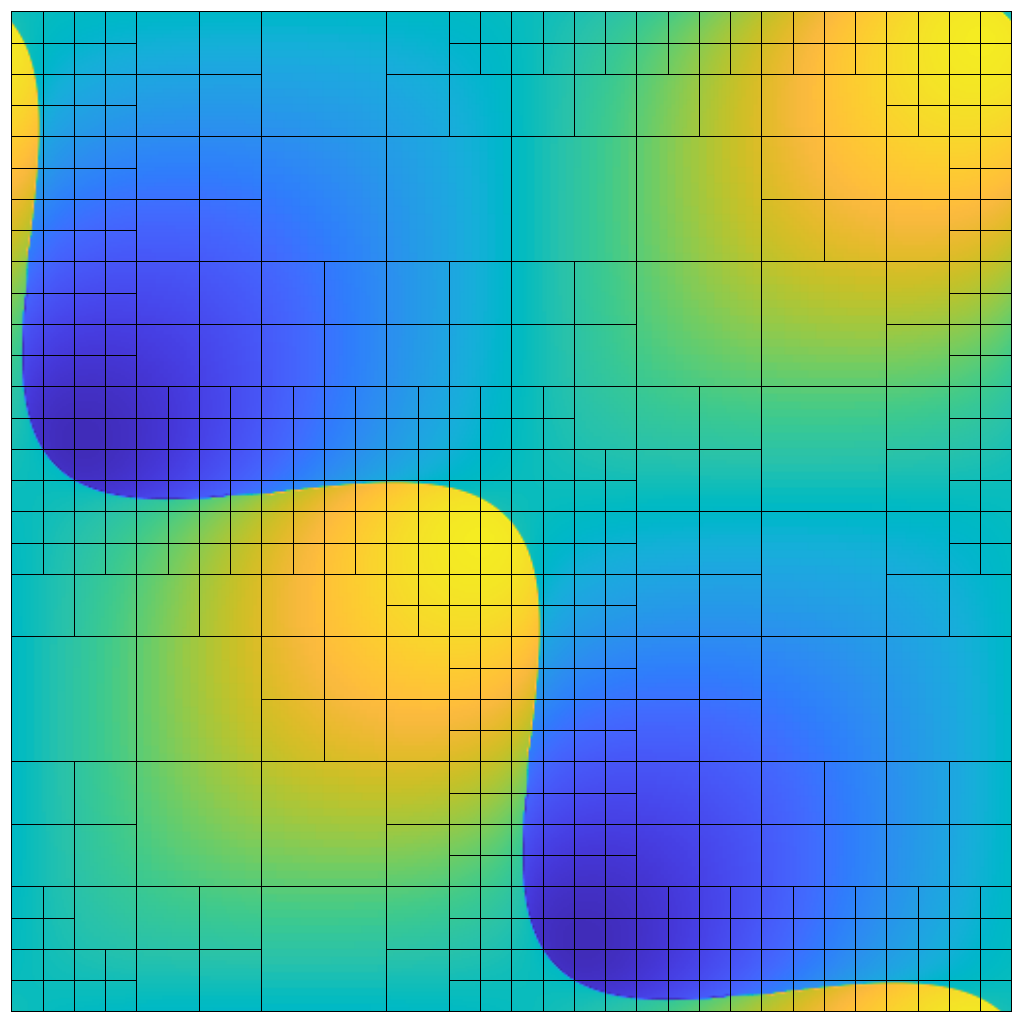} \hfil
\includegraphics[width=0.28\textwidth]{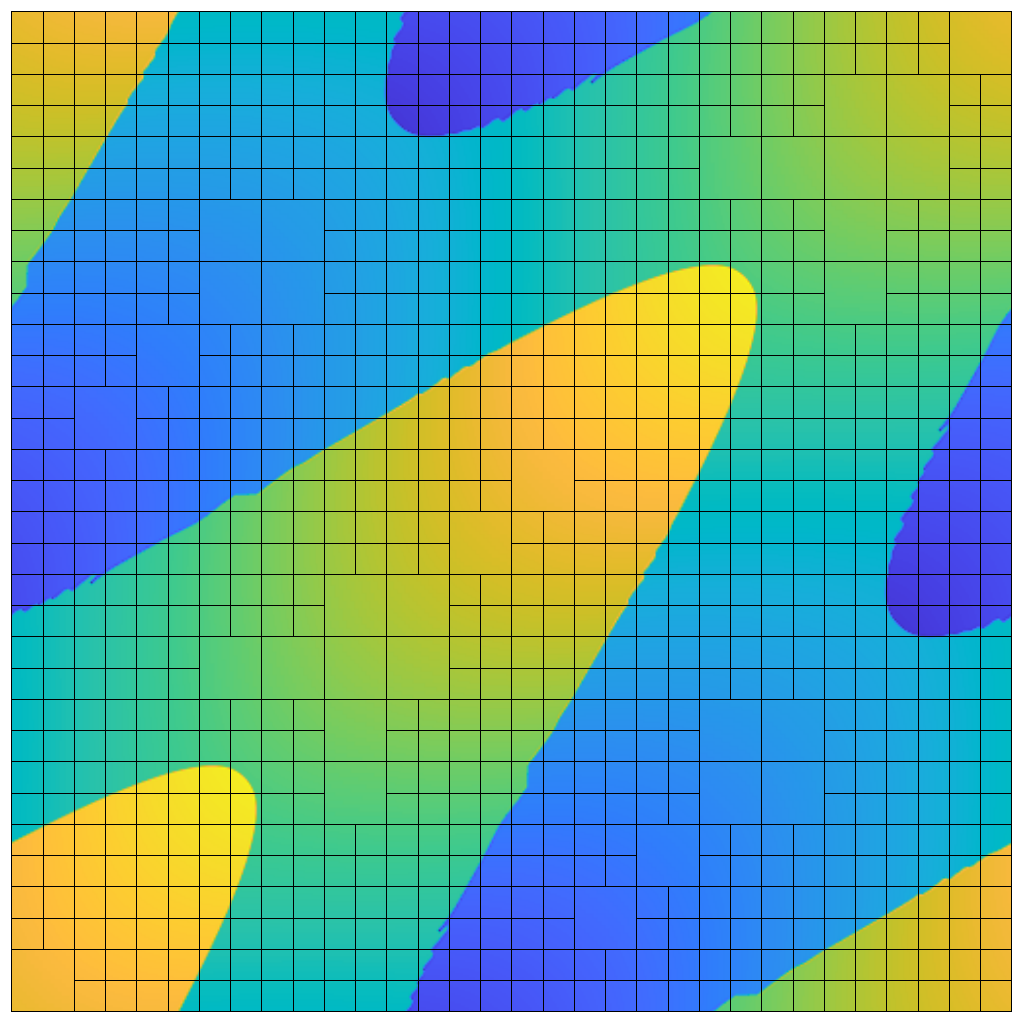}\hfil
\includegraphics[width=0.28\textwidth]{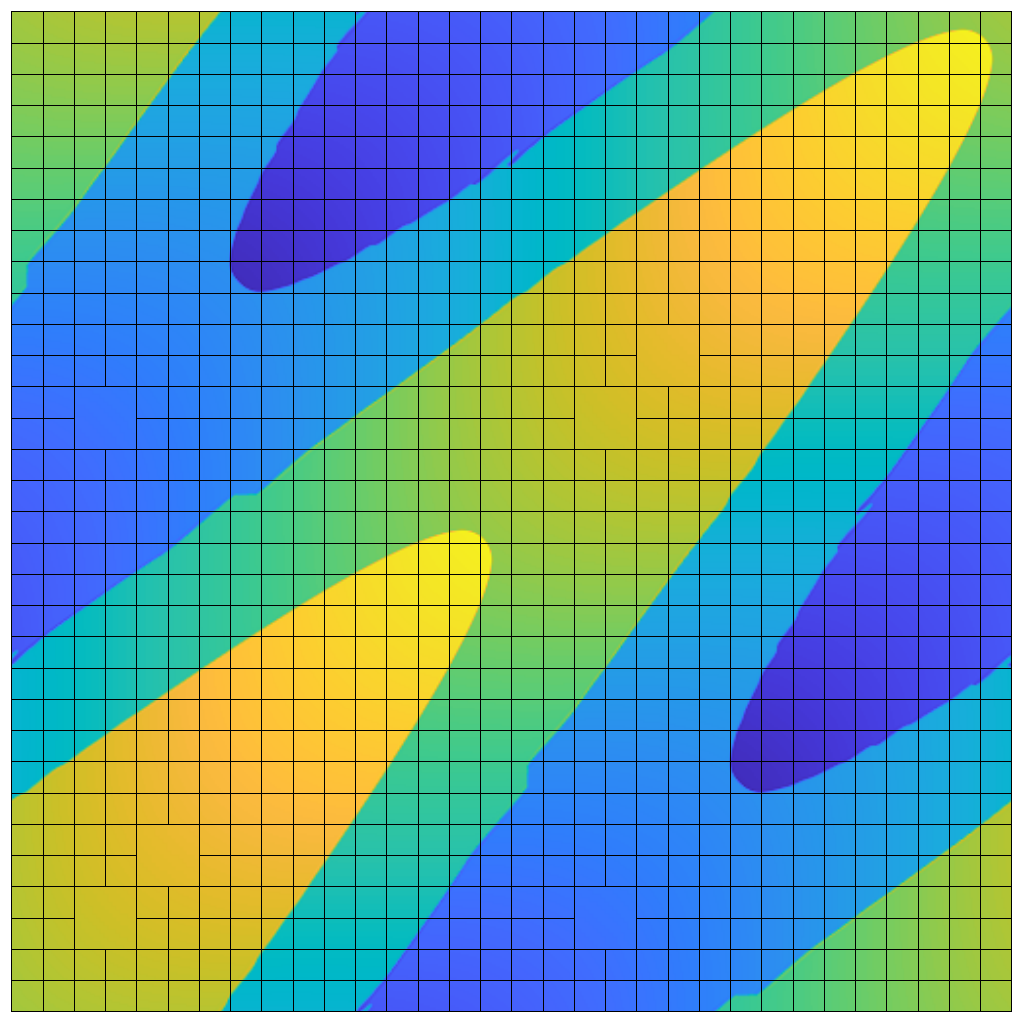} \\
\vspace{0.5cm}
\includegraphics[width=0.28\textwidth]{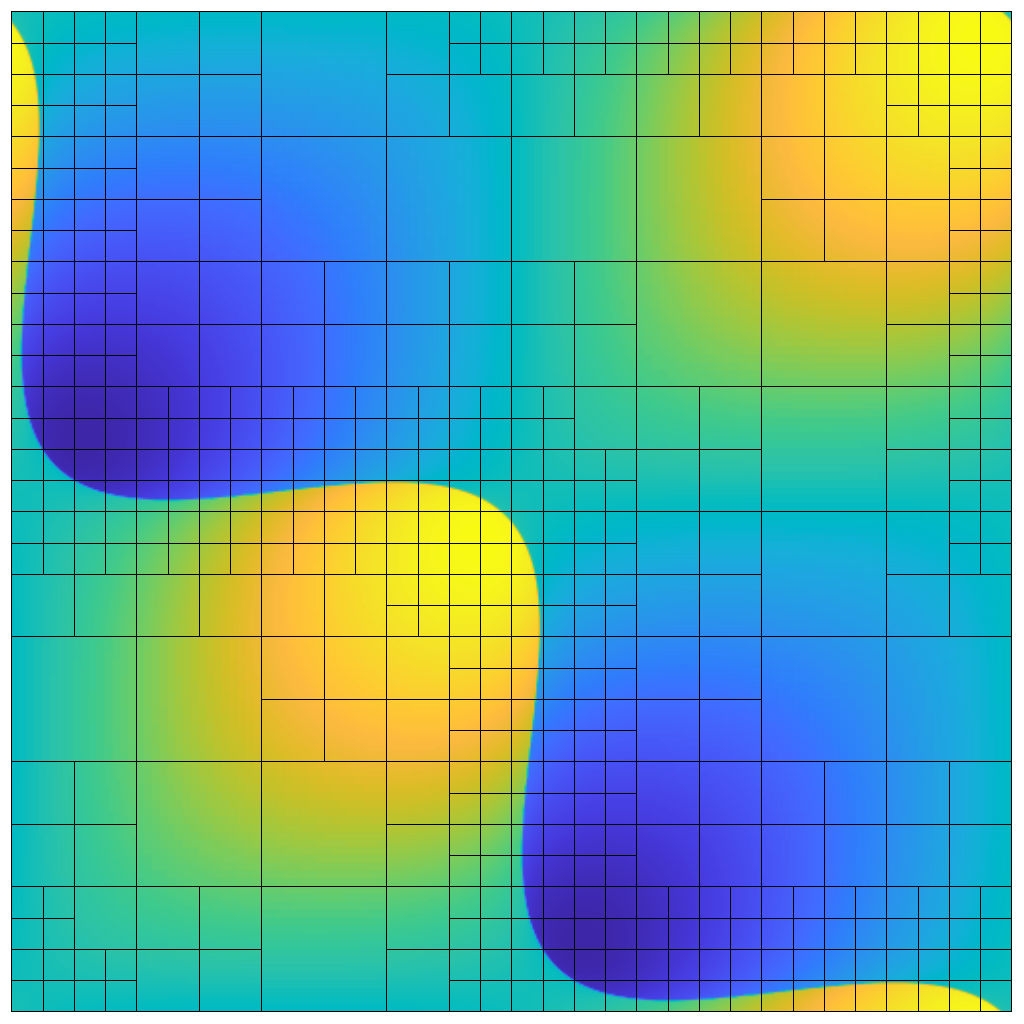} \hfil
\includegraphics[width=0.28\textwidth]{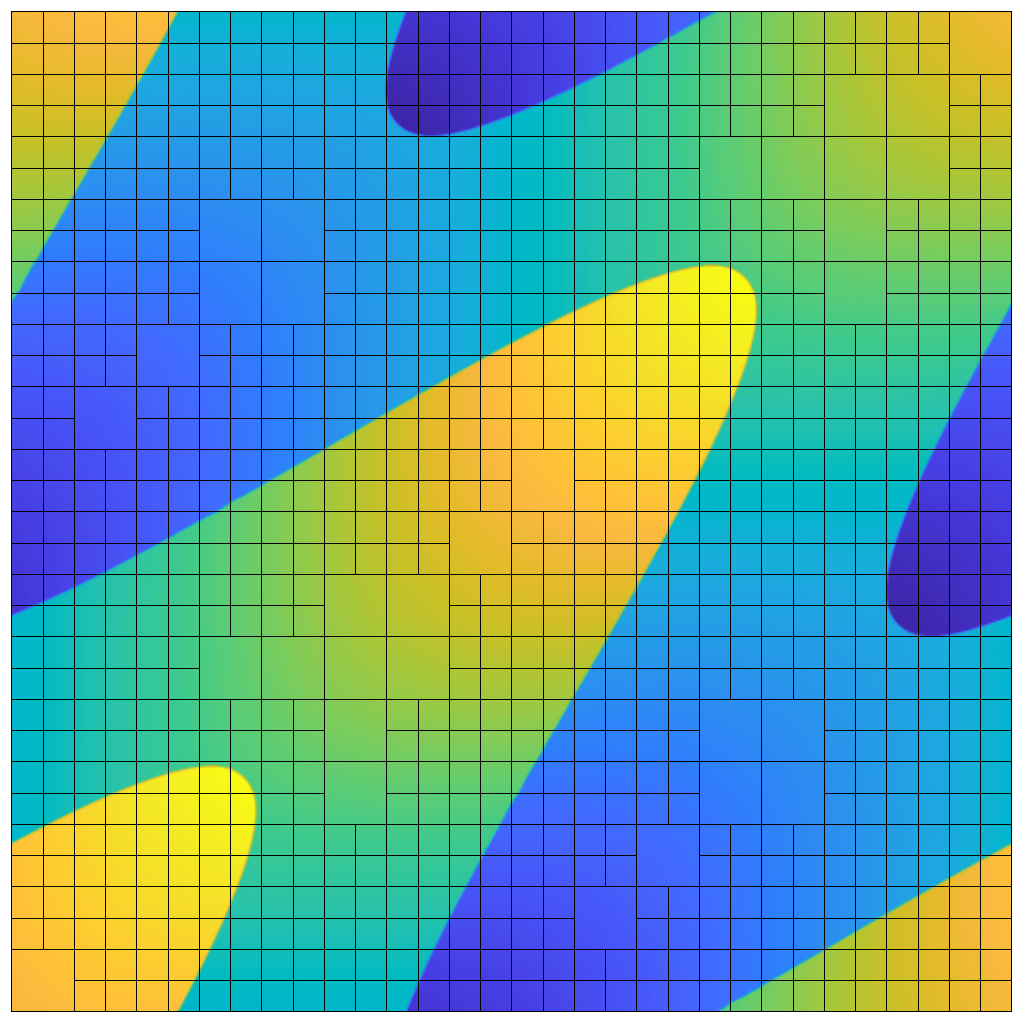}\hfil
\includegraphics[width=0.28\textwidth]{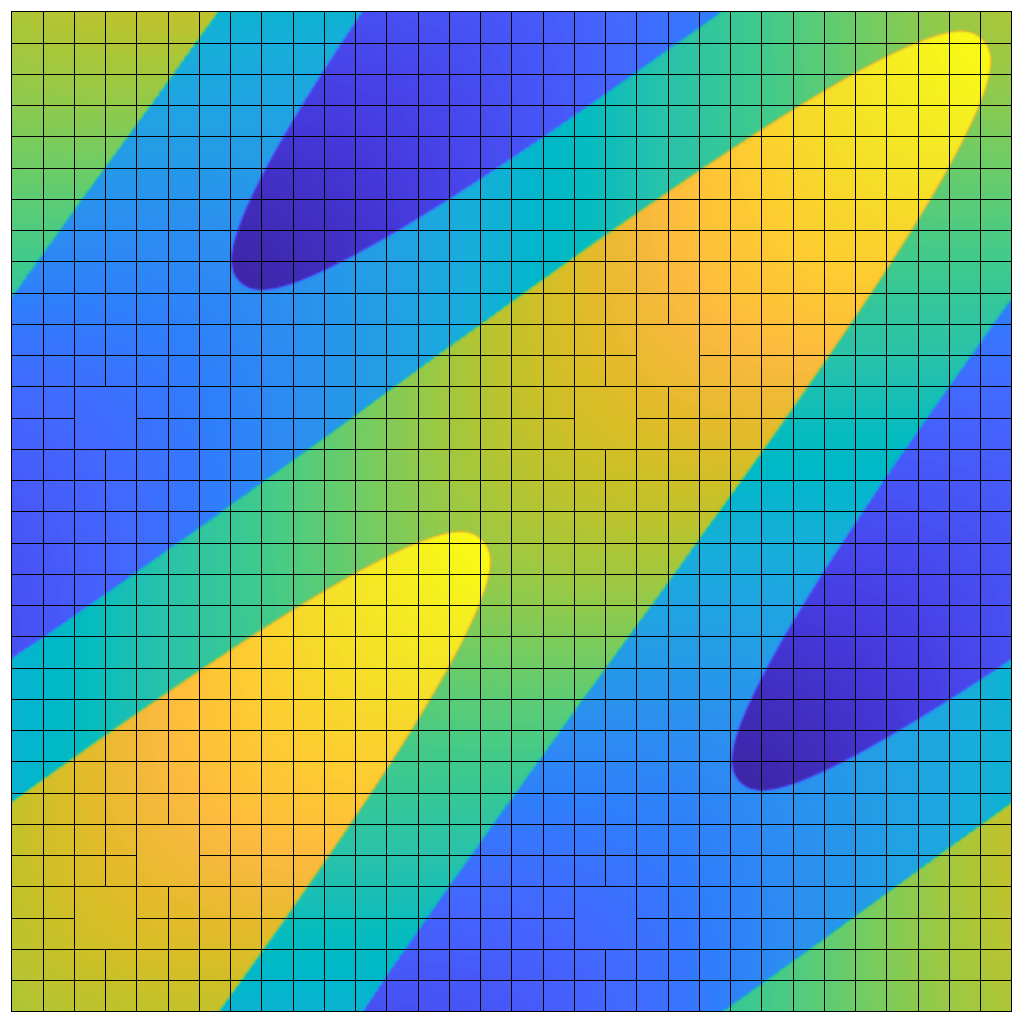}
\end{center}
	\caption{Solution to (\ref{eqn:sinsin01}) at times
		$t = 0.2$ (left),
		$t = 0.6$ (middle) and 
		$t = 1.0$ (right) without limiter (top) and with limiter (bottom) for refinement level $3-5$.} \label{Burgers}
\end{figure}

At later times shocks arise and we therefore also used the bound preserving limiter introduced in \cite{article:CHK2021}. Note that the characteristic speed changes sign which has been observed to lead to some numerical difficulties as explained in \cite{article:HKS2019,article:CHK2021}. Our unlimited method shows some unphysical oscillations  along the curve where the characteristic speeds changes sign as can be seen in Figure \ref{Burgers} (top). The use of the bound preserving limiter avoids these inaccuracies. In this simulation grids on levels $3-5$ are used.

As a refinement criteria we search for steep gradients, so a patch is refined, if $$\frac{q_{i+1,j}-q_{i-1,j}}{2\Delta x} \ge 15 \text{ \ or \ } \frac{q_{i,j+1}-q_{i,j+1}}{2\Delta y} \ge 15 \ ,$$ for any $i,j$.  We also use the refinement criteria to determine whether we need to apply limiting.

\subsection{Acoustics}
\subsubsection{Convergence study for acoustics}
In order to investigate the order of convergence of the adaptive Active Flux method for acoustics, we consider a test problem from Luk\'{a}c\v{o}v\'{a} et al. \cite{article:LMW2000}, for  which the exact solution is explicitly known. In this test problem the acoustic equations (\ref{eqn:acoustics}) with initial values of the form
\begin{equation}
  \begin{aligned}
p(x,y,0) & = -\frac{1}{c} \left( (\sin(2 \pi x) + \sin(2 \pi y) \right)
\\
u(x,y,0) & = 0 \\
v(x,y,0) & = 0
  \end{aligned}
\end{equation} are considered 
on the domain $[-1,1] \times [-1,1]$. The speed of sound is set to $c=1$ and periodic boundary conditions are imposed. The exact solution has the form
\begin{equation*}
  \begin{split}
p(x,y,t) & = - \frac{1}{c} \cos(2 \pi c t) \left( \sin(2 \pi x) +
  \sin(2 \pi y) \right) \\
u(x,y,t) & = \frac{1}{c} \sin(2 \pi c t) \cos(2 \pi x) \\
v(x,y,t) &= \frac{1}{c} \sin(2 \pi c t) \cos(2 \pi y). 
  \end{split}
\end{equation*}
We compute numerical solutions at time $t=1$ using time steps, which satisfy the condition $\mbox{CFL} \le 0.4$. For this smooth solution, adaptive mesh refinement does not really make sense. However,  in order to test the accuracy of the adaptive method, we require refinement along the diagonal of the domain independently of the solution structure. Each patch is a Cartesian grid with $16 \times 16$ grid cells.
\begin{figure}[t]
\begin{center}
\includegraphics[width=0.35\textwidth]{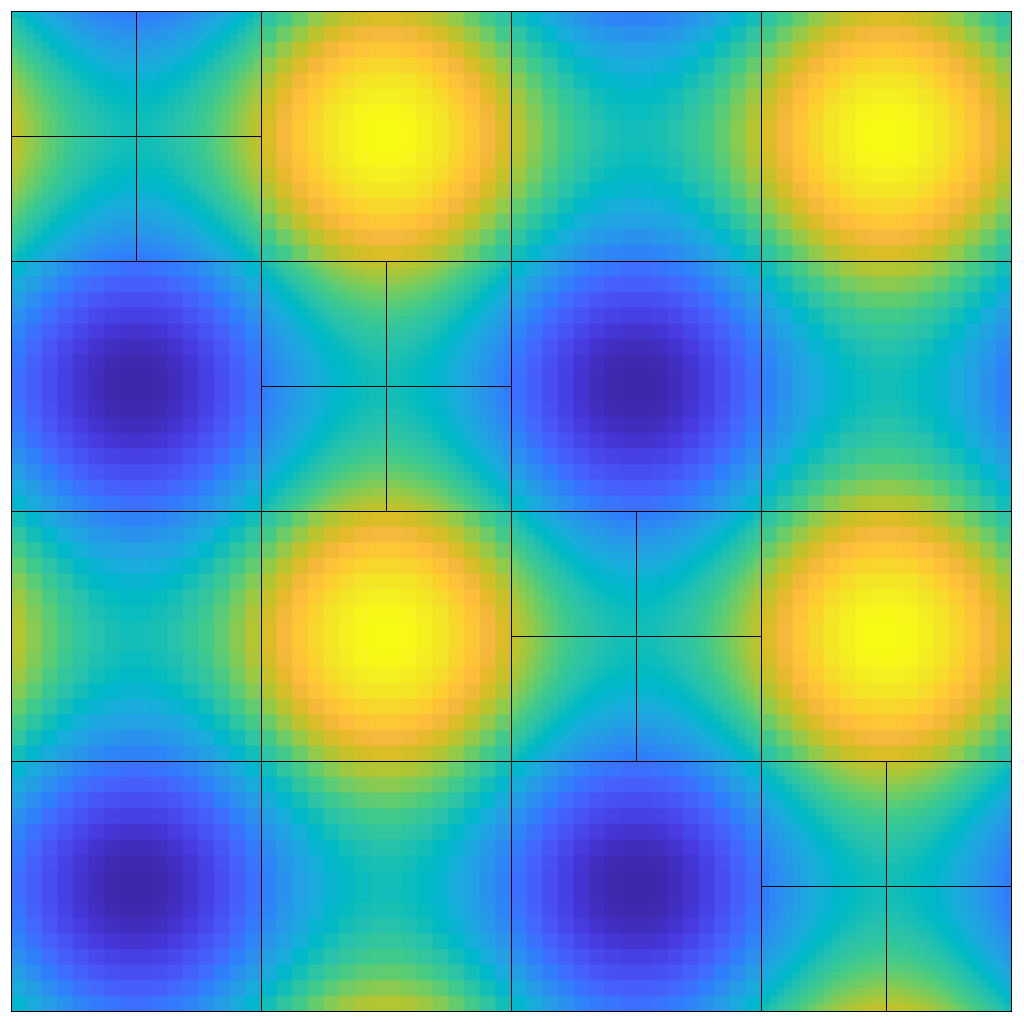}\hfil
\includegraphics[width=0.35\textwidth]{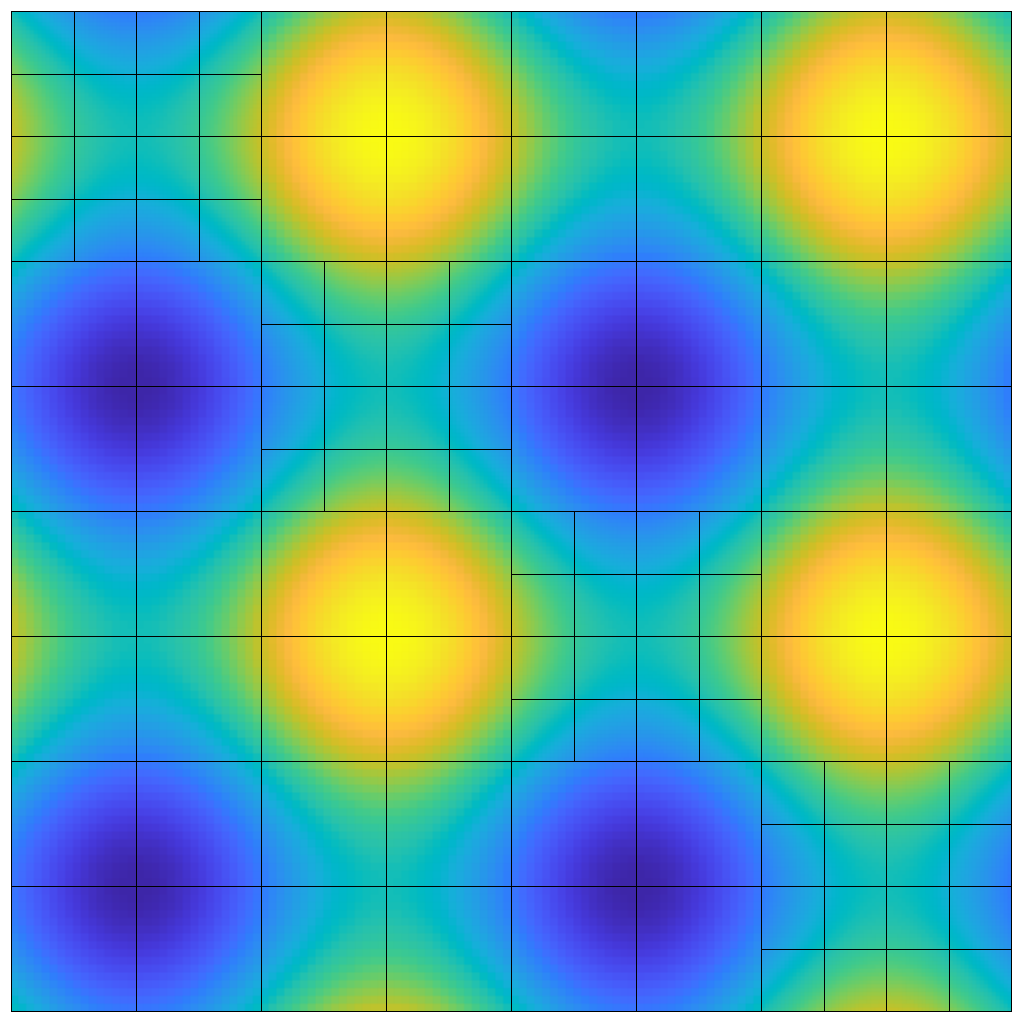}
\end{center}
\caption{Pressure and grid patches for the acoustics problem from
  Luk\'{a}c\v{o}v\'{a} et al. \cite{article:LMW2000}
  with refinement levels 2-3 (left)
  and 3-4 (right) at $t=1.0$.}
\end{figure}
We compare grids with refinement level $2-3$  and $3-4$ with results obtained on a uniform grid with level $2$ and $3$. The results of our convergence study  for pressure and velocity are shown in Table \ref{Tab:acousticEOC}.

\begin{table}[htb]
    \caption{Error at time $t=1$ measured in $\| \cdot \|_1$-norm  and EOC for the Luk\'{a}c\v{o}v\'{a} test problem.}
\sisetup{
scientific-notation = true, 
round-mode=places, round-precision=2}
\begin{center}
\begin{tabular}{
    *1{S[table-column-width=1.5cm,table-text-alignment=center]}
    *2{S[table-column-width=2.25cm,table-text-alignment=center]}
    *2{S[table-column-width=1.25cm,table-text-alignment=center, round-precision=4]}
}
\toprule
{Level}   & \multicolumn{2}{c}{Error}   & \multicolumn{2}{c}{EOC}  \\ 
\midrule
          & {$p$}   & {$u,v$}           & {$p$}  & {$u,v$} \\ 
\midrule
{2} & \num{2.231166e-04} & \num{2.520880e-05} &  {---}   & {---}\\
{3} & \num{2.793825e-05} & \num{2.648027e-06} & \num{2.9975} &  \num{3.2509}\\ 
\midrule
{2-3} & \num{1.641979e-04} & \num{3.328607e-05} & {---}    & {---} \\
{3-4} & \num{2.057597e-05} & \num{3.869926e-06} & \num{3.1045} & \num{2.8564}\\ 
\bottomrule                       
  \end{tabular}  
  \label{Tab:acousticEOC}
\end{center}
\end{table}

The accuracy obtained on the adaptively refined mesh is comparable with the accuracy obtained on a regular grid which uses the coarser grid in the whole domain. This test shows that the the accuracy is maintained at the interface between the coarse and the fine grids. 

\subsubsection{High frequency acoustics}
\begin{figure}[htb]
\begin{center}
\includegraphics[width=0.35\textwidth]{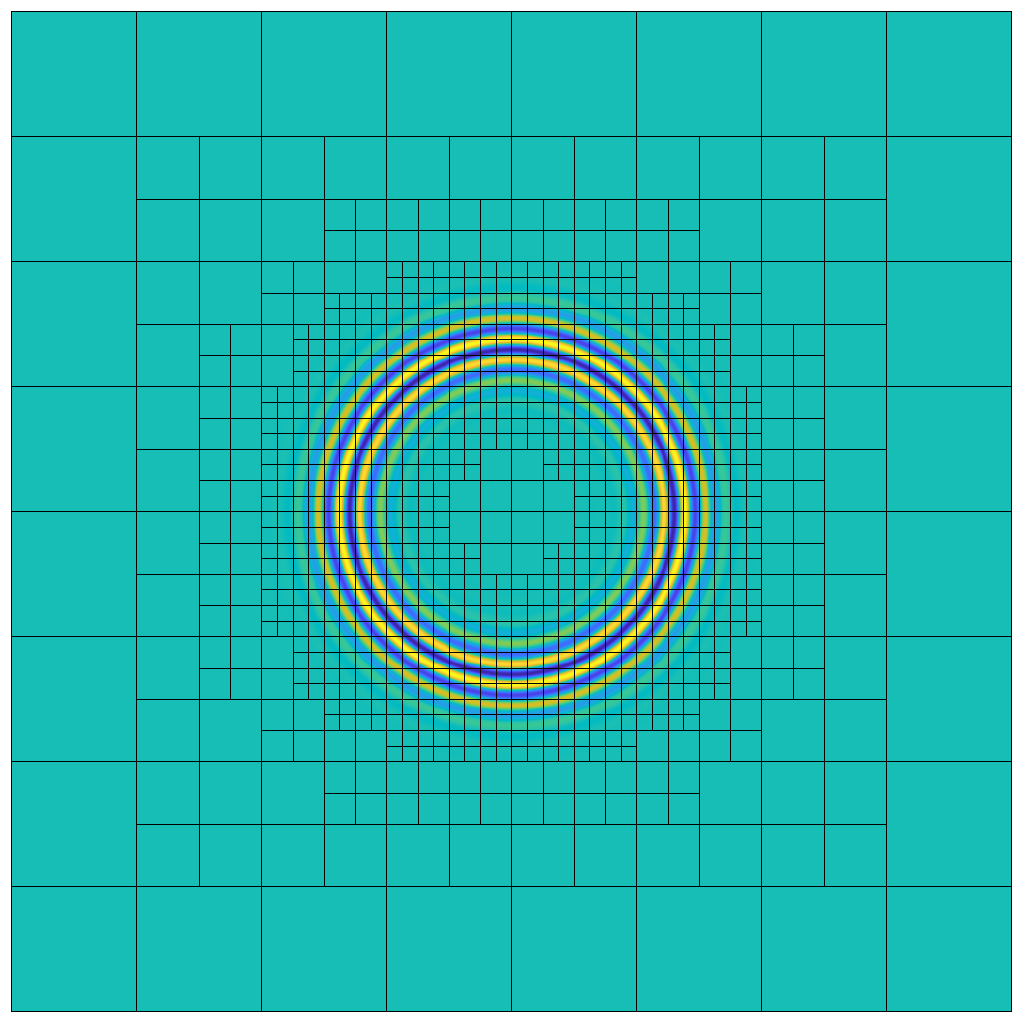} \hfil
\includegraphics[width=0.35\textwidth]{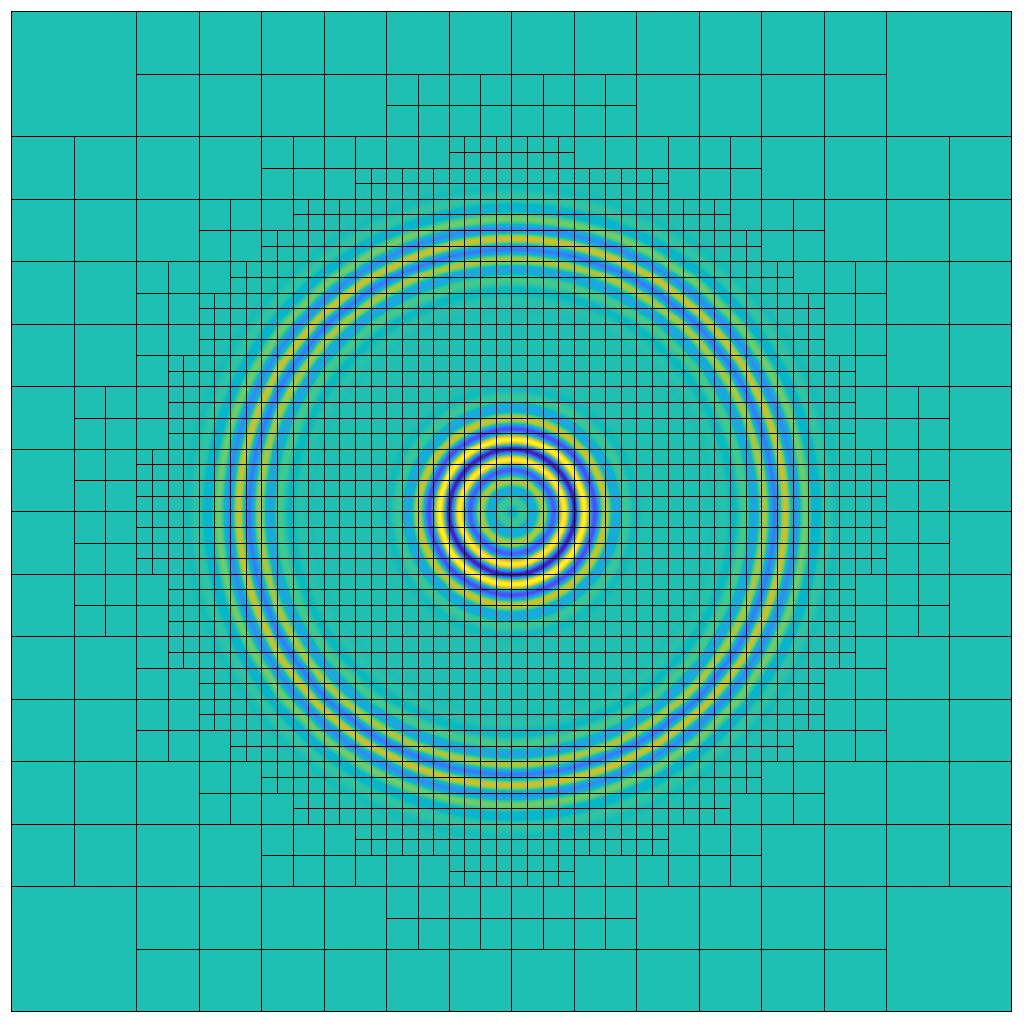} \\
\vspace{0.5cm}
\includegraphics[width=0.35\textwidth]{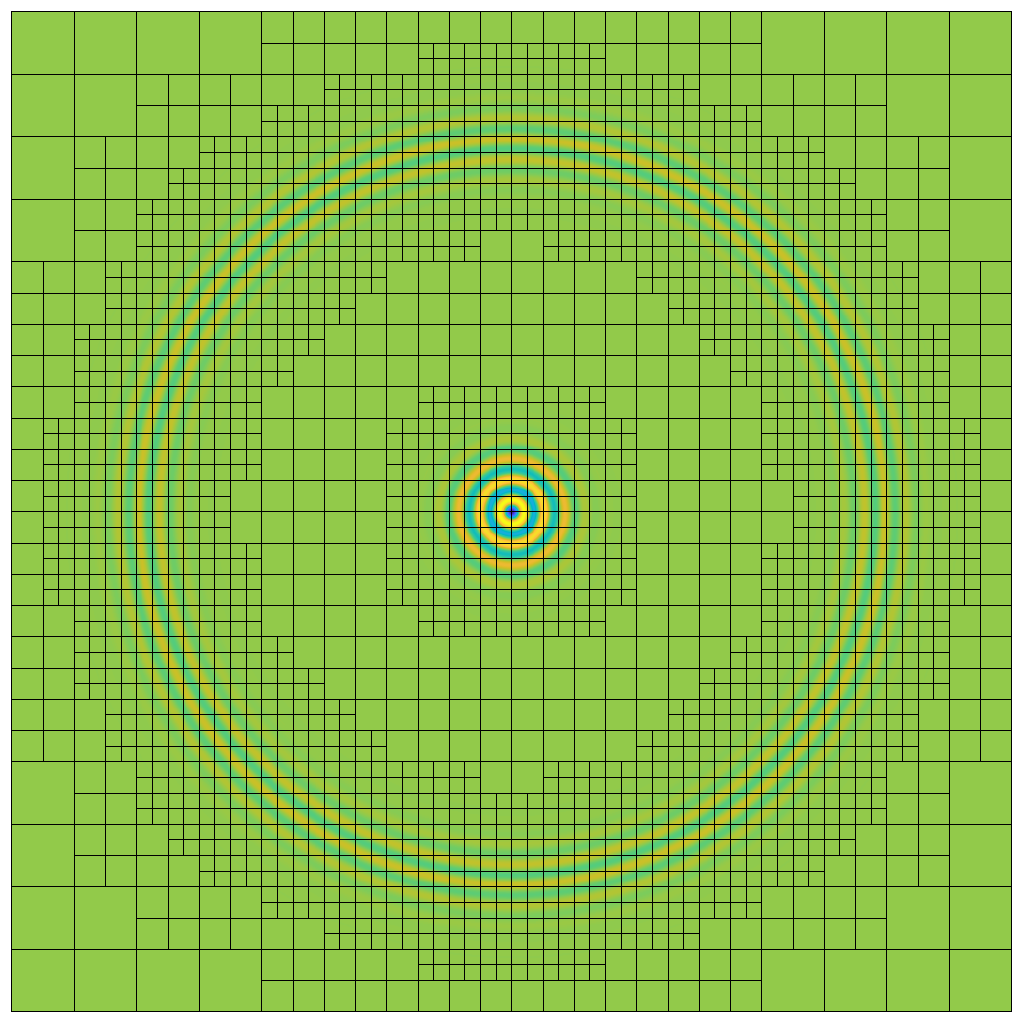} \hfil
\includegraphics[width=0.35\textwidth]{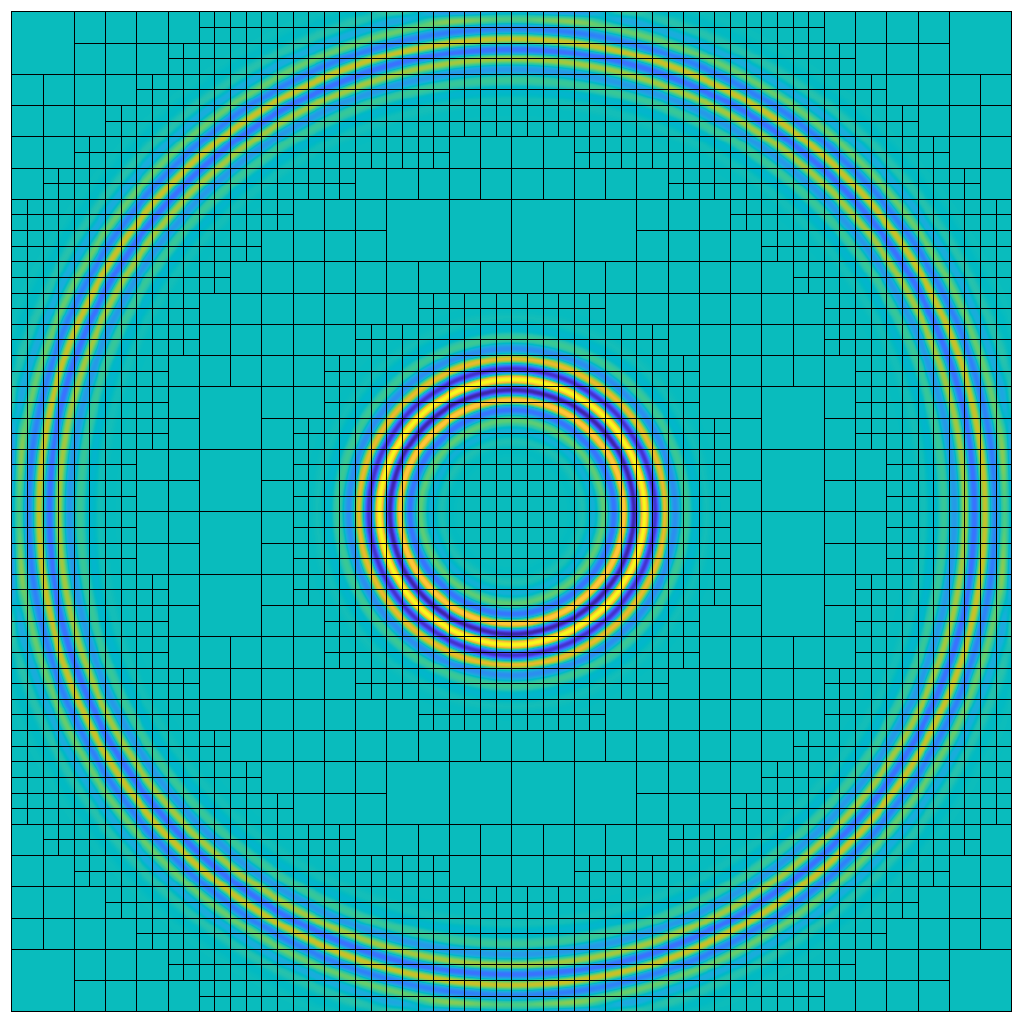}
\end{center}
	\caption{Plots of pressure of the high frequency acoustic wave at times
		$t = 0.0$ (top left),
		$t = 0.3$ (top right),  
		$t = 0.6$ (bottom left) and $t = 0.9$ (bottom
                right). In the bottom left plot the pressure at the
                  center is much larger and therefore a different
                  color map is used to visualize the solution structure.
                } 
    \label{HighFreq}
  \end{figure}

\begin{figure}[htb]
  \begin{center}
\includegraphics[width=0.485\textwidth,clip=true,trim=3cm 0cm 2cm 0cm]{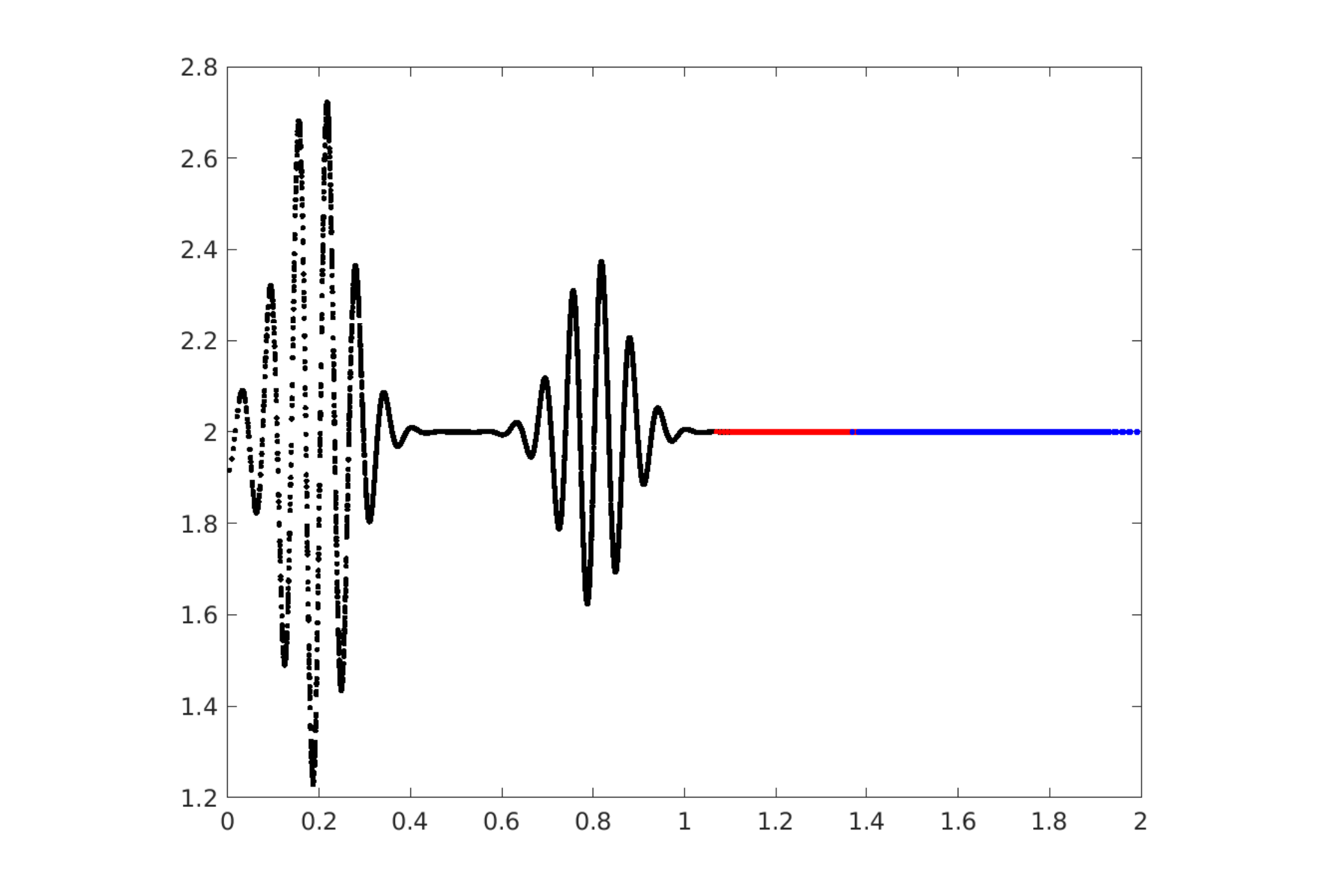} \hfil
\includegraphics[width=0.485\textwidth,clip=true,trim=3cm 0cm 2cm 0cm]
{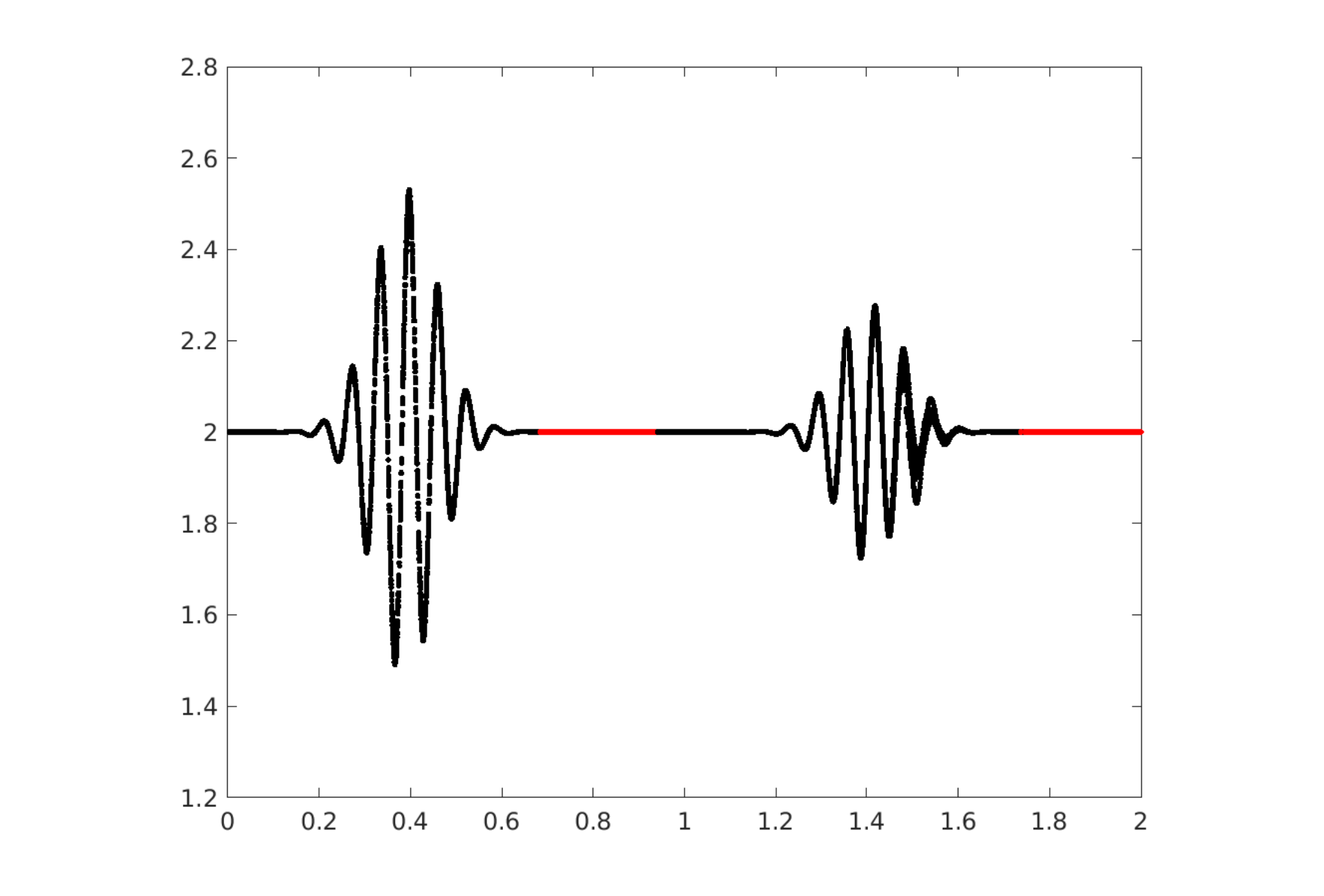} 
\end{center}
\caption{Scatter plots of pressure for the high frequency acoustic
  wave at time $t=0.3$ and $t=0.9$. }
\label{fig:HighFreq-scatter}
\end{figure}

Now we consider a test problem where adaptive mesh refinement allows a more efficient computation of the solution structure.  We consider the acoustic equations (\ref{eqn:acoustics}) with initial condition
\begin{equation}\label{eqn:highfreq}
p(x,y,0) = 2+\exp(-100(r(x,y)-0.5)^2))\sin(100r(x,y)), \ \ u(x,y,0)=0, \ \ v(x,y,0)=0
\end{equation} 
on the domain $[-1.5,1.5]\times[-1.5,1.5]$.  Initially a circular shaped acoustic wave with high frequency pressure oscillations is given. This leads to acoustic waves moving outwards and inwards. The inward moving acoustic wave gets reflected in the center of the domain (around the time shown in the third plot) and afterwards propagates outwards. Adaptive mesh refinement is used to resolve these high frequency waves. Figure \ref{HighFreq} shows solutions at times $t=0,0.3,0.6,0.9$ with refinement level $3-6$. A patch is refined if
 $|p_{max}-p_{min}| \ge 0.001$.
In Figure \ref{fig:HighFreq-scatter} shows scatter plots of the solution at two different times. There are no visable spurious grid effects apart from some smearing of the solution structure caused by the boundary of the computational domain. For this simulation periodic boundary conditions were used.

\subsection{Approximation of steady states}
Finally, we consider the approximation of a stationary vortex as described by Barsukow et al.\ \cite{article:Barsukow2020,article:BHKR2018}. The initial values have the form
\begin{equation*}
  \begin{split}
    p(r,0) & = 0 \\
    \vec{u}(r) & = \vec{n} \left\{ \begin{array}{ccc}
                                       5r & : & 0 \le r \le 0.2 \\
                                       2-5r & : & 0.2 < r \le 0.4 \\
                                       0 & : & r>0.4,\end{array}\right.
  \end{split}
\end{equation*}
with $r = \sqrt{x^2 + y^2}$, $\vec{n} = (-\sin \phi, \cos \phi)^T$, $\phi \in [0, 2\pi)$ and $\vec{u} = (u,v)^T$. In order to test how well the method preserves the steady state, we compute the numerical solution at time $t=100$. The AMR computation uses grids of level 3-5 with  $16 \times 16$ grid cells on each patch.  

\begin{figure}[htb]
\begin{center}
\includegraphics[width=0.35\textwidth]{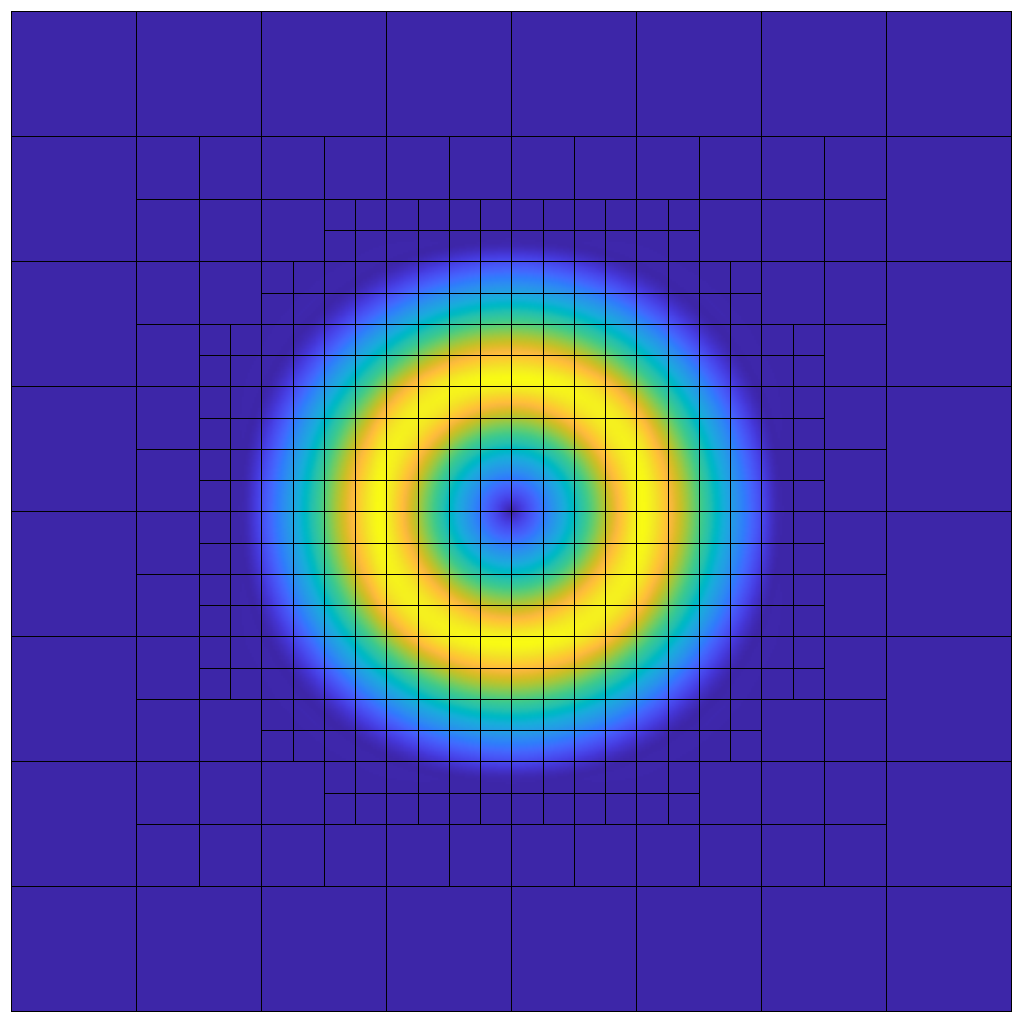} \hfil
\includegraphics[height=0.25\textheight]{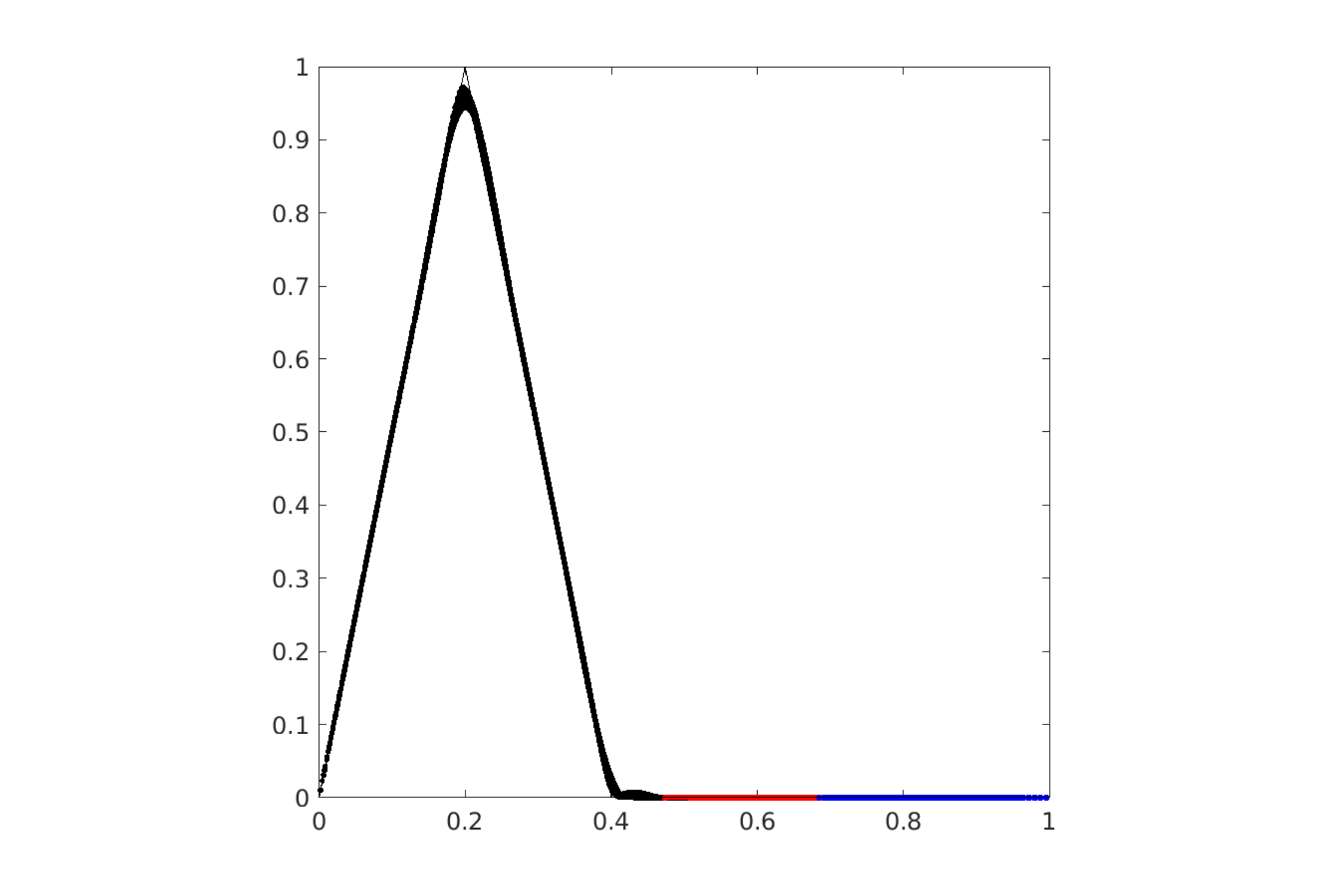}
\end{center}
\caption{\label{fig:steadystate} Computation of the stationary vortex
  at time $t=100$ using an adaptively refined grid of level 3-5. The
  left plot shows $|\vec{u}|$, the right plot shows a scatter plot of $|\vec{u}|$.}
\end{figure}

In \cite{article:Barsukow2020,article:BHKR2018}, the authors showed that the Cartesian grid Active Flux method, with the evolution operator presented in \cite{article:BHKR2018}, is stationary preserving. The method used here does not have this property. Nevertheless, the third order accuracy of the Active Flux method together with adaptive mesh refinement leads to accurate approximations as shown in Figure \ref{fig:steadystate}. 
We also tested the evolution operator from \cite{article:BHKR2018}. While we could confirm the preservation of steady states on a regular Cartesian grid, using the solver with AMR leads to an instability.

\section*{Conclusions}
We showed that the Active flux method can be used on Cartesian grids with adaptive mesh refinement and subcycling. The transfer of grid information between different Cartesian grid patches can be implemented without loss of third order accuracy by making use of the degrees of freedom of the Active Flux method. Our approach benefits from the local stencil of the Active Flux method.

For advective transport, new Active Flux methods have been presented which preserve constant states on regular Cartesian grids and grids with adaptive mesh refinement without subcycling. In practical computations with subcycling we also observed good accuracy although constant states are not exactly preserved.

The AMR  concept of the Active Flux method can also be used for two-dimensional linear hyperbolic systems as illustrated for the acoustic equations. However, the preservation of steady states, a property that was recently shown for the Active Flux method on regular Cartesian grids \cite{article:Barsukow2020}, does not carry over to adaptively refined meshes with the approaches presented in this paper.

\section*{Data availability statement}
For our implementations we used ForestClaw, which is publicly available on GitHub. We added the Active Flux method as a new solver in ForestClaw. 


\bibliographystyle{plain}

\end{document}